\documentclass[10pt]{amsart}

\usepackage{amssymb}
\usepackage{amsmath,amscd}
\usepackage{amsthm}
\usepackage{graphicx,psfrag}

\theoremstyle{plain}
\newtheorem{theorem}{Theorem}[section]

\newtheorem{lemma}[theorem]{Lemma}
\newtheorem{corollary}[theorem]{Corollary}
\newtheorem{proposition}[theorem]{Proposition}

\theoremstyle{definition}
\newtheorem{definition}[theorem]{Definition}

\theoremstyle{remark}
\newtheorem{remark}[theorem]{Remark}
\newtheorem*{notation}{Notation}
\newtheorem*{acknowledgments}{Acknowledgments}

\newtheorem{operation}{Operation}

\numberwithin{figure}{section}

\begin{document}

\title{Heegaard surfaces and measured laminations, I: the Waldhausen conjecture}

\author{Tao Li}
\thanks{Partially supported by NSF grants DMS-0102316 and DMS-0406038} 

\address{Department of Mathematics \\
 Oklahoma State University \\
 Stillwater, OK 74078\\
 USA}
\curraddr{Department of Mathematics \\
 Boston College \\
 Chestnut Hill, MA 02467\\
 USA}
\email{taoli@bc.edu}
%\url{www.math.okstate.edu/~tli}

\begin{abstract}
We give a proof of the so-called generalized Waldhausen conjecture, which says that an orientable irreducible atoroidal 3--manifold has only finitely many Heegaard splittings in each genus, up to isotopy.  Jaco and Rubinstein have announced a proof of this conjecture using different methods.  
\end{abstract}

\maketitle

\begin{psfrags}

\setcounter{tocdepth}{1}
\tableofcontents

\section{Introduction}
A Heegaard splitting of a closed orientable 3--manifold $M$ is a decomposition of $M$ into two handlebodies along an embedded surface called a Heegaard surface.  Heegaard splittings were introduced to construct and classify 3--manifolds. Every 3--manifold has a Heegaard splitting, and one can construct a Heegaard splitting of arbitrarily large genus for any 3--manifold by adding trivial handles to a Heegaard splitting.  An important problem in 3--manifold topology is the classification of Heegaard splittings of a 3--manifold.  The main questions are whether there are different Heegaard splittings in a 3--manifold and how the different Heegaard splittings are related.  A conjecture of Waldhausen asserts that a closed orientable 3--manifold has only a finite number of Heegaard splittings of any given genus, up to homeomorphism.  Johannson \cite{Jo1,Jo2} proved this conjecture for Haken manifolds.  If $M$ contains an incompressible torus, one may construct an infinite family of homeomorphic but non-isotopic Heegaard splittings using Dehn twists along the torus.  The so-called generalized Waldhausen conjecture says that a closed, orientable and atoroidal 3--manifold has only finitely many Heegaard splittings of any fixed genus, up to isotopy.  This is also proved to be true for Haken manifolds by Johannson \cite{Jo1,Jo2}.  The main purpose of this paper is to prove the generalized Waldhausen conjecture.

\begin{theorem}\label{main}
A closed, orientable, irreducible and atoroidal 3--manifold has only finitely many Heegaard splittings in each genus, up to isotopy.
\end{theorem}

Jaco and Rubinstein have announced a proof using normal surface theory and 1--efficient triangulations.  The main tools used in this paper are measured laminations and branched surfaces.  In a sequel to this paper \cite{L4}, we use measured laminations and Theorem~\ref{T2} of this paper to prove a much stronger result for non-Haken 3--manifolds, which says that, for non-Haken manifolds, adding trivial handles is virtually the only way of creating new Heegaard splittings.  

Methods of laminations and branched surfaces have been very useful in solving some seemingly unrelated problems, such as \cite{L5, L2}.  This is the first time that they are used on Heegaard splittings.  Both \cite{L4} and this paper use branched surfaces to analyze Heegaard surfaces.  The main technical issues in this paper are on measured laminations, whereas the arguments in \cite{L4} rely more on the properties of strongly irreducible Heegaard splittings.

A theorem of Schleimer \cite{Sch} says that every Heegaard splitting of sufficiently large genus has the disjoint curve property.  So, an immediate corollary is that $M$ contains only finitely many full Heegaard splittings, see \cite{Sch}.

\begin{corollary}
In any closed orientable 3--manifold, there are only finitely many full Heegaard splittings, up to isotopy.
\end{corollary}

Theorem~\ref{main} provides a well-known approach to understand the structure of the mapping class group of 3--manifolds.  Conjecturally, the mapping class group for such a 3--manifold is finite, but it is not clear how to obtain a geometric description of the mapping class group.  For instance, there is no example of a non-trivial element in the mapping class group of such a 3--manifold that is invariant on a strongly irreducible Heegaard splitting.  Very recently, Namazi \cite{Na} used the result of this paper and showed that if the distance of a Heegaard splitting is large then the mapping class group is finite.

We give a very brief outline of the proof.  By a theorem of Rubinstein and Stocking \cite{St}, every strongly irreducible Heegaard splitting is isotopic to an almost normal surface.  So, similar to \cite{FO}, one can construct a finite collection of branched surfaces using normal disks and almost normal pieces of a triangulation, such that every almost normal strongly irreducible Heegaard surface is fully carried by a branched surface in this collection.  If no branched surface in this collection carries any surface with non-negative Euler characteristics, then Theorem~\ref{main} follows immediately from a simple argument of Haken in normal surface theory.  The key of the proof is to show that one can split a branched surface into a finite collection of branched surfaces so that no branched surface in this collection carries any normal torus and up to isotopy, each almost normal Heegaard surface is still carried by a branched surface in this collection, see sections \ref{Storus} and \ref{Slam}.  Most of the paper are dedicated to proving Theorem~\ref{T2}, and Theorem~\ref{main} follows easily from this theorem, see section~\ref{Sproof}.

\begin{theorem}\label{T2}
Let $M$ be a closed orientable, irreducible and atoroidal 3--manifold, and suppose $M$ is not a small Seifert fiber space.  Then, $M$ has a finite collection of branched surfaces, such that
\begin{enumerate}
\item each branched surface in this collection is obtained by gluing together normal disks and at most one almost normal piece, similar to \cite{FO},
\item up to isotopy, each strongly irreducible Heegaard surface is fully carried by a branched surface in this collection,
\item no branched surface in this collection carries any normal 2--sphere or normal torus. 
\end{enumerate} 
\end{theorem} 

In the proof, we also use some properties of 0--efficient triangulations \cite{JR}.  The use of 0--efficient triangulations does not seem to be absolutely necessary, but it makes many arguments much simpler.  Jaco and Rubinstein also have a theory of 1--efficient triangulations, which can simplify our proof further, but due to the status of their paper, we decide not to use it.  Some arguments in this paper are also similar in spirit to those in \cite{AL, L1}.  One can also easily adapt the arguments in this paper into \cite{AL} so that the algorithm in \cite{AL} works without the use of 1--efficient triangulations.

\begin{acknowledgments}
I would like to thank Bus Jaco for many conversations and email communications on their theory of efficient triangulations.  I also thank Saul Schleimer and Ian Agol for helpful conversations. I would also like to thank the referee for many corrections and suggestions.
\end{acknowledgments}

\section{Heegaard splittings, almost normal surfaces and branched surfaces}\label{Spre}

\begin{notation}
Throughout this paper, we will denote the interior of $X$ by $int(X)$, the closure (under path metric) of $X$ by $\overline{X}$, and the number of components of $X$ by $|X|$.  We will use $\eta(X)$ to denote the closure of a regular neighborhood of $X$.  We will use $M$ to denote a closed, orientable, irreducible and atoroidal 3--manifold, and we always assume $M$ is not a Seifert fiber space.
\end{notation}

In this section, we will explain some basic relations between Heegaard splittings, normal surface theory and branched surfaces.  We will also explain some terminology and operations that are used throughout this paper.

\subsection{Heegaard splittings}

A \emph{handlebody} is a compact 3--manifold homeomorphic to a regular neighborhood of a connected graph embedded in $\mathbb{R}^3$.  A \emph{Heegaard} splitting of a closed 3--manifold $M$ is a decomposition $M=H_1\cup_SH_2$, where $S=\partial H_1=\partial H_2=H_1\cap H_2$ is a closed embedded separating surface and each $H_i$ ($i=1,2$) is a handlebody.  The surface $S$ is called a \emph{Heegaard surface}, and the genus of $S$ is the genus of this Heegaard splitting.  The boundary of a regular neighborhood of the 1--skeleton of any triangulation of $M$ is a Heegaard surface.  Hence, any closed orientable 3--manifold has a Heegaard splitting.  The notion of Heegaard splitting can be generalized to manifolds with boundary, but in this paper, we only consider Heegaard splittings of closed 3--manifolds.

Heegaard splitting became extremely useful when Casson and Gordon introduced strongly irreducible Heegaard splitting.

\begin{definition}
A \emph{compressing disk} of a handlebody $H$ is a properly embedded disk in $H$ with boundary an essential curve in $\partial H$.  A Heegaard splitting is \emph{reducible} if there is an essential curve in the Heegaard surface that bounds compressing disks in both handlebodies.  A Heegaard splitting $M=H_1\cup_SH_2$ is  \emph{weakly reducible} \cite{CG} if there exist a pair of compressing disks $D_1\subset H_1$ and $D_2\subset H_2$ such that $\partial D_1\cap\partial D_2=\emptyset$.  If a Heegaard splitting is not reducible (resp. weakly reducible), then it is \emph{irreducible} (resp. \emph{strongly irreducible}).
\end{definition}

A closed 3--manifold $M$ is \emph{reducible} if $M$ contains an embedded 2--sphere that does not bound a 3--ball.  A lemma of Haken \cite{H} says that if $M$ is reducible, then every Heegaard splitting is reducible.  Casson and Gordon \cite{CG} showed that if a Heegaard splitting of a non-Haken 3--manifold is irreducible, then it is strongly irreducible.   

The following theorem of Scharlemann \cite{S} is useful in proving Theorem~\ref{main}.

\begin{theorem}[Theorem 3.3 of \cite{S}]\label{Tsch}
Suppose $H_1\cup_S H_2$ is a strongly irreducible Heegaard splitting of a 3--manifold and $V\subset M$ is a solid torus such that $\partial V$ intersects $S$ in parallel essential non-meridian curves.  Then $S$ intersects $V$ in a collection of $\partial$--parallel annuli and possibly one other component, obtained from one or two $\partial$-parallel annuli by attaching a tube along an arc parallel to a subarc of $\partial V$.
\end{theorem}

\subsection{Almost normal surfaces}

A normal disk in a tetrahedron is either a triangle cutting off a vertex or a quadrilateral separating two opposite edges, see Figure 3 of \cite{JR} for a picture. An almost normal piece in a tetrahedron is either an octagon, or an annulus obtained by connecting two normal disks using an unknotted tube, see Figures 1 and 2 in \cite{St} for pictures.

\begin{definition}
Suppose a 3--manifold $M$ has a triangulation $\mathcal{T}$.  We use $\mathcal{T}^{(i)}$ to denote the $i$--skeleton of $\mathcal{T}$.  Let $S$ be a surface in $M$ that does not meet the 0--skeleton $\mathcal{T}^{(0)}$ and is transverse to $\mathcal{T}^{(1)}$ and $\mathcal{T}^{(2)}$.  $S$ is called a \emph{normal surface} (or we say $S$ is normal) with respect to $\mathcal{T}$ if the 2--skeleton $\mathcal{T}^{(2)}$ cuts $S$ into a union of normal disks.  $S$ is called an \emph{almost normal surface} if $S$ is normal except in one tetrahedron $T$, where $T\cap S$ consists of normal disks and at most one almost normal piece.
\end{definition}

Rubinstein and Stocking \cite{R,St} (see also \cite{K}) showed that any strongly irreducible Heegaard surface is isotopic to an almost normal surface with respect to any triangulation of the 3--manifold. 

Normal surfaces, introduced by Kneser \cite{Kn}, have been very useful in the study of incompressible surfaces.  The results and techniques in normal surface theory are similarly applicable to almost normal surfaces.  

Let $S$ be a surface in $M$ transverse to the 1--skeleton of  $\mathcal{T}$ and with $S\cap\mathcal{T}^{(0)}=\emptyset$.  We define the \emph{weight} (or the combinatorial area) of $S$, denoted by $weight(S)$, to be $|S\cap\mathcal{T}^{(1)}|$.  Let $\alpha$ be an arc such that $\alpha\cap\mathcal{T}^{(1)}=\emptyset$ and $\alpha$ is transverse to $\mathcal{T}^{(2)}$.  We define the combinatorial length of $\alpha$, denoted by $length(\alpha)$, to be $|\alpha\cap\mathcal{T}^{(2)}|$.  After a small perturbation, we may assume any arc to be disjoint from $\mathcal{T}^{(1)}$ and transverse to $\mathcal{T}^{(2)}$.  In this paper, when we mention the length of an arc, we always use such combinatorial length.

Let $S$ be a closed embedded normal surface in $M$. If we cut $M$ open along $S$, the manifold with boundary $\overline{M-S}$ has an induced cell decomposition.  One can also naturally define normal disks and normal surfaces in $\overline{M-S}$ with respect to this cell decomposition.  An embedded disk in a 3--cell is a \emph{normal disk} if its boundary curve does not meet the 0--cells, meets at least one edge, and meets no edge more than once.  

An isotopy of $M$ is called a \emph{normal isotopy} if it is invariant on the cells, faces, edges and vertices of the triangulation.  In this paper, we will consider two normal surfaces (or laminations) the same if they are isotopic via a normal isotopy.  Up to normal isotopy there are only finitely many equivalence classes of normal disks, and these are called \emph{normal disk types}.  There are 7 types of normal disks in a tetrahedron.

\subsection{Branched surfaces}

A \emph{branched surface} in $M$ is a union of finitely many compact smooth surfaces glued together to form a compact subspace (of $M$)  locally modeled on Figure~\ref{branch}(a). 

Given a branched surface $B$ embedded in a 3--manifold $M$, we denote by $N(B)$ a regular neighborhood of $B$, as shown in Figure~\ref{branch}(b).  One can regard $N(B)$ as an $I$--bundle over $B$, where $I$ denotes the interval $[0,1]$.  We denote by $\pi : N(B)\to B$ the projection that collapses every $I$--fiber to a point.  The \emph{branch locus} of $B$ is $L=\{b\in B:$ $b$ does not have a neighborhood homeomorphic to $\mathbb{R}^2  \}$.  So, $L$ can be considered as a union of smoothly immersed curves in $B$, and we call a point in $L$ a \emph{double point} of $L$ if any small neighborhood of this point is modeled on Figure~\ref{branch}(a).  We call the closure (under the path metric) of each component of $B-L$ a \emph{branch sector} of $B$.  We say that a surface (or lamination) $S$ is carried by a branched surface $B$ (or carried by $N(B)$) if $S$ lies in $N(B)$ and is transverse to the $I$--fibers of $N(B)$.  We say $S$ is \emph{fully carried} by $B$, if $S\subset N(B)$ transversely intersects every $I$--fiber of $N(B)$.  The boundary of $N(B)$ consists of two parts, the horizontal boundary, denoted by $\partial_hN(B)$, and the vertical boundary, denoted by $\partial_vN(B)$.  The vertical boundary is a union of subarcs of $I$--fibers of $N(B)$ and the horizonal boundary is transverse to the $I$--fibers of $N(B)$, as shown in Figure~\ref{branch} (b).

\begin{figure}
\begin{center}
\psfrag{(a)}{(a)}
\psfrag{(b)}{(b)}
\psfrag{horizontal}{$\partial_hN(B)$}
\psfrag{v}{$\partial_vN(B)$}
\includegraphics[width=4.0in]{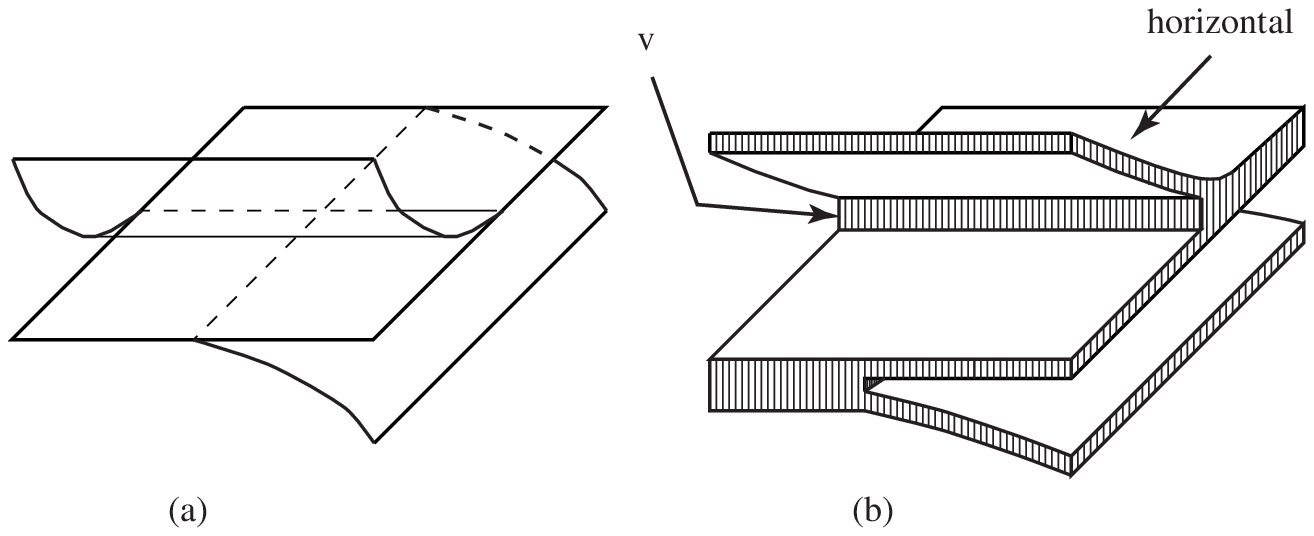}
\caption{}\label{branch}
\end{center}
\end{figure}

Let $\mu\subset N(B)$ be a lamination carried by $N(B)$ (or $B$), and let $b$ be a branch sector of $B$.  We say that $\mu$ \emph{passes through} the branch sector $b$ if $\mu\cap\pi^{-1}(int(b))\ne\emptyset$, where $\pi : N(B)\to B$ is the collapsing map.  So, $\mu$ is fully carried by $B$ if and only if $\mu$ passes through every branch sector.  Let $x\in int(b)$ be a point and $I_x=\pi^{-1}(x)$  the corresponding $I$--fiber.  If $\mu$ is a closed surface, then $m=|I_x\cap\mu|$ is a non-negative integer and $m$ does not depend on the choice of $x\in int(b)$. We call $m$ the weight of $\mu$ at the branch sector $b$.

\begin{definition}
A \emph{disk of contact} is an embedded disk in $N(B)$ transverse to the $I$--fibers of $N(B)$ and with $\partial D\subset\partial_vN(B)$, see \cite{FO} for a picture.  A \emph{monogon} is a disk $E$ properly embedded in $M-int(N(B))$ with $\partial E=\alpha\cup\beta$, where $\alpha\subset\partial_vN(B)$ is a subarc of an $I$--fiber of $N(B)$ and $\beta\subset\partial_hN(B)$.  If a component of $M-int(N(B))$ is a 3--ball whose boundary consists of two disk components of $\partial_hN(B)$ and a component of $\partial_vN(B)$, then we call this 3--ball a $D^2\times I$ region.  If a component of $M-int(N(B))$ is a solid torus, whose boundary consists of an annulus component of $\partial_hN(B)$ and a component of $\partial_vN(B)$, and a meridian disk of the solid torus is a monogon, then we call this solid torus a $monogon\times S^1$ region.  Let $A$ be an annulus in $N(B)$.  We call $A$ a \emph{vertical annulus} if $A$ is a union of subarcs of the $I$--fibers of $N(B)$.
\end{definition}

For any embedded (almost) normal surface $S$, by identifying all the normal disks of the same disk type as in \cite{FO}, we obtain a branched surface fully carrying $S$.  Since $M$ is compact and there are only finitely many different types of normal disks, there are only finitely many such branched surfaces.   This construction is first used by Floyd and Oertel \cite{FO} to study incompressible surfaces, then used in \cite{G8, L1, L2, L3} to study essential laminations and immersed surfaces.  Since strongly irreducible Heegaard surfaces are isotopic to almost normal surfaces, by the argument above, we have the following.

\begin{proposition}\label{Pfinite}
There is a finite collection of branched surfaces in $M$ with the following properties.
\begin{enumerate}
\item each branched surface is obtained by gluing normal disks and at most one almost normal piece, similar to \cite{FO},
\item after isotopy, every strongly irreducible Heegaard surface is fully carried by a branched surface in this collection.
\end{enumerate}
\end{proposition}\qed

Let $B$ be a branched surface, and $\mathcal{B}$ the set of branched sectors of $B$.  If a subset of $\mathcal{B}$ also form a branched surface $B'$, then we call $B'$ a \emph{sub-branched surface} of $B$.  If a lamination $\mu$ is carried but not fully carried by $B$, then the branch sectors that $\mu$ passes through form a sub-branched surface of $B$ that fully carries $\mu$.

Let $B$ be a branched surface as in Proposition~\ref{Pfinite}.  If $B$ does not contain any almost normal piece, then every surface carried by $B$ is a normal surface.  Suppose $B$ contains an almost normal branched sector, which we denote by $b_A$.  Let $S$ be an almost normal surface fully carried by $B$.  By the definition of almost normal surface, the weight of $S$ at the branch sector $b_A$ is one.  Therefore, it is easy to see that $B_N=B-int(b_A)$ is a sub-branched surface of $B$.  We call $B_N$ the \emph{normal part} of $B$.  Every surface carried by $B_N$ is a normal surface.

In this paper, we assume all the 3--manifolds are orientable.  So, if $S$ is a non-orientable surface carried by $N(B)$, then a small neighborhood of $S$ in $N(B)$ is a twisted $I$--bundle over $S$ and the boundary of this twisted $I$--bundle is an orientable surface carried by $B$.  Thus, we have the following trivial proposition.

\begin{proposition}\label{Pnonori}
If a branched surface in an orientable 3--manifold does not carry any 2--sphere (resp. torus), then $B$ does not carry any projective plane (resp. Klein bottle).
\end{proposition}

\subsection{Splitting branched surfaces}\label{SSsplit}

\begin{definition}
An isotopy of $N(B)$ is called a $B$--\emph{isotopy} if it is invariant on each $I$--fiber of $N(B)$.  We say two surfaces carried by $N(B)$ are $B$--isotopic if they are isotopic via a $B$--isotopy of $N(B)$. 
\end{definition}

Let $S$ be a compact surface embedded in $N(B)$ transverse to the $I$--fibers of $N(B)$, and let $N_S$ be a closed neighborhood of $S$ in $N(B)$.  We call $N_S$ a \emph{fibered neighborhood} of $S$ in $N(B)$ if $N_S$ is an $I$--bundle over $S$ with each $I$--fiber of $N_S$ a subarc of an $I$--fiber of $N(B)$.  After some small perturbation, $N(B)-int(N_S)$ can be considered as a fibered neighborhood $N(B')$ of another branched surface $B'$.  We say that $B'$ is obtained by \emph{splitting $B$ along $S$}.  For most splittings considered in this paper, we have $\partial S\cap\partial_vN(B)\ne\emptyset$ and $S$ is orientable.  If $\mu$ be a surface or lamination carried by $N(B)$ and $S\subset N(B)-\mu$ is $B$--isotopic to a sub-surface of (a leaf of) $\mu$, then we also say that $B'$ is obtained by splitting $B$ along $\mu$.   The inverse operation of splitting is called \emph{pinching}, and we say that $B$ is obtained by pinching $B'$.  If $B'$ is a branched surface obtained by splitting $B$, then we may naturally consider $N(B')$ as a subset of $N(B)$ with the induced fiber structure.  For any lamination $\mu$ carried by $B$, we say that $\mu$ is also carried by $B'$ if after some $B$--isotopies, $\mu$ is carried by $N(B')$ with $\mu\subset N(B')\subset N(B)$.

Suppose $B'$ is obtained by splitting $B$.  Since we can regard $N(B')\subset N(B)$, we have the following obvious proposition.

\begin{proposition}\label{Pob}
Suppose $B'$ is obtained by splitting $B$.  Then, any lamination carried by $B'$ is also carried by $B$.
\end{proposition}\qed

The converse of Proposition~\ref{Pob} is not true.  It is possible that some lamination is carried by $B$ but not carried by $B'$.  For example, in Figure~\ref{sp1}, the train track $\tau_2$ is obtained by splitting the train track $\tau$ on the top.  However, any lamination fully carried by $\tau_1$ or $\tau_3$ is carried by $\tau$ but not carried by $\tau_2$.  Nevertheless, every lamination carried by $\tau$ is carried by some $\tau_i$ ($i=1,2,3$).  Moreover, $\tau_2$ is a sub-traintrack of each $\tau_i$.  

One can apply such different splittings (as in Figure~\ref{sp1}) to branched surfaces.  The next proposition is also obvious, see section~\ref{Storus} for a more general discussion of such splittings.

\begin{proposition}\label{Psplit}
Let $B$ be a branched surface and $\{S_n\}$ a sequence of distinct closed surfaces fully carried by $B$.  Suppose $B'$ is a branched surface obtained by splitting $B$ and $B'$ fully carries some $S_m$.  Then, there is a finite collection of branched surfaces, such that 
\begin{enumerate}
\item each branched surface in this collection is obtained by splitting $B$, and $B'$ is in this collection,
\item each $S_n$ is fully carried by a branched surface in this collection,
\item if another branched surface $B''$ in this collection carries $S_m$, then $B'$ is a sub-branched surface of $B''$.
\end{enumerate}
\end{proposition}
\begin{proof} 
First note that, in the one-dimension lower example Figure~\ref{sp1}, if $B$ is $\tau$ and $B'$ is $\tau_i$ ($i=1,2,3$), then $\tau_1$, $\tau_2$ and $\tau_3$ form a collection of train tracks satisfying the 3 conditions of the proposition.  The 2-dimensional case is similar.  Any splitting can be viewed as a sequence of successive local splittings similar to Figure~\ref{sp1}.  During each local splitting, one can enumerate all possible splittings as in Figure~\ref{sp1} and get a collection of branched surfaces satisfying the conditions of this proposition.
\end{proof}
\begin{remark}
If $B$ is obtained by gluing normal disks and at most one almost normal piece as in Proposition~\ref{Pfinite}, then the branched surface after splitting is also obtained by gluing normal disks and almost normal pieces.  Moreover, if $\{S_n\}$ in Proposition~\ref{Psplit} are almost normal surfaces, since each $S_n$ has at most one almost normal piece, we may assume that each branched surface in this collection has at most one branch sector containing an almost normal piece.
\end{remark}

\begin{figure}
\begin{center}
\psfrag{sp}{splitting}
\psfrag{t}{$\tau$}
\psfrag{1}{$\tau_1$}
\psfrag{2}{$\tau_2$}
\psfrag{3}{$\tau_3$}
%\psfrag{4}{$\tau_4$}
\includegraphics[width=4in]{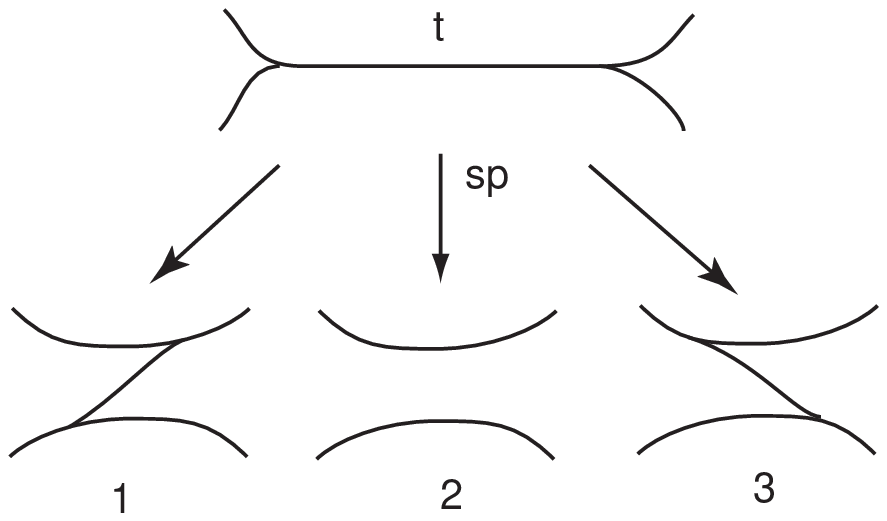}
\caption{} \label{sp1}
\end{center}
\end{figure}

\section{Measured laminations and projective lamination spaces}\label{SPL}

Let $B$ be a branched surface in $M$, and $F\subset N(B)$ be a surface carried by $B$.  Let $L$ be the branch locus of $B$, and suppose $b_1,\dots, b_N$ are the components of $B-L$.  For each $b_i$, let $x_i=|F\cap\pi^{-1}(b_i)|$.  One can describe $F$ using a non-negative integer point $(x_1,\dots,x_N)\in\mathbb{R}^N$, and  $(x_1,\dots,x_N)$ is a solution to the system of branch equations of $B$, see \cite{FO,O} for more details.  $F$ is fully carried by $B$ if and only if each $x_i$ is positive. Each branch equation is of the form $x_k=x_i+x_j$.  We use $\mathcal{S}(B)$ to denote the set of non-negative solutions to the system of branch equations of $B$.  This gives a one-to-one correspondence between closed surfaces carried by $B$ and integer points in $\mathcal{S}(B)$. Throughout this paper, we do not distinguish a surface carried by $B$ from its corresponding non-negative integer point $(x_1,\dots,x_N)\in\mathcal{S}(B)$.  We will call $x_n$ the \emph{weight} (or the \emph{coordinate}) of the surface at the branch sector corresponding to $b_n$.

Let $F_1$ and $F_2$ be embedded closed orientable surfaces carried by $N(B)$ and suppose $F_1\cap F_2\ne\emptyset$.  In general, there are two directions to perform  cutting and pasting along an intersection curve of $F_1\cap F_2$, but only one of them results in surfaces still transverse to the $I$--fibers of $N(B)$.  We call such cutting and pasting the \emph{canonical} cutting and pasting.  This is similar to the Haken sum in normal surface theory.  We use $F_1+F_2$ to denote the surface after the canonical cutting and pasting.  This is a very natural operation, because if $F_1=(x_1,\dots,x_N)$ and $F_2=(y_1,\dots,y_N)$ in $\mathcal{S}(B)$ then $F_1+F_2=(x_1+y_1,\dots,x_N+y_N)$.  Moreover, this sum preserves the Euler characteristic, $\chi(F_1)+\chi(F_2)=\chi(F_1+F_2)$.

A theorem of Haken \cite{Ha} says that there is a finite set of fundamental integer solutions $F_1,\dots,F_k$ in $\mathcal{S}(B)$, such that any integer solution in $\mathcal{S}(B)$ can be written as $\sum_{i=1}^kn_iF_i$, where each $n_i$ is a non-negative integer.  In other words, every surface carried by $B$ can be obtained by the canonical cutting and pasting on multiple copies of $F_1,\dots,F_k$.  So, if $B$ does not carry any 2--sphere or torus, by Proposition~\ref{Pnonori}, $B$ does not carry any surface with non-negative Euler characteristic and hence there are only finitely many surfaces (carried by $B$) with any given genus.

The positive non-integer points of $\mathcal{S}(B)$ correspond to measured laminations fully carried by $B$.  We refer to \cite{H, O, MS} for details.  Roughly speaking, one can construct the measured lamination as follows, see \cite{O} and section 2 of \cite{H}.  We can first pinch each component of $\partial_vN(B)$ to a circle and change $N(B)$ to $N_w(B)$, see Figure 1.2 of \cite{O} or Figure 2.2 of \cite{L1}. $N(B)$ is basically the same as $N_w(B)$ except the vertical boundary of $N(B)$ becomes the cusp of $N_w(B)$.  For each branch sector of $B$, we can take an $I$--bundle over this sector with a standard horizontal foliation.  For any positive point in $\mathcal{S}(B)$, when we glue the branch sectors together, we glue the foliations according to the weights at these sectors, see Figure 1.2 of \cite{O}.  This produces a singular foliation of $N_w(B)$ where the cusps are the singularity.  So, there are a finite number of singular leaves.  Now, one can split $B$ along these singular leaves.  This is usually an infinite process, and the inverse limit is a measured lamination fully carried by $B$.  It is not hard to show that if the singular foliation does not contain any compact leaf, then the singular leaves are dense in the lamination (see \cite{H, MS}).  Throughout this paper, we always assume our measured laminations are constructed in this fashion.  So, we may assume that there is a one-to-one correspondence between a point in $\mathcal{S}(B)$ and a measured lamination carried by $B$.

Measured laminations in 3--manifolds have many remarkable properties.  We say a lamination is \emph{minimal} if it has no sub-lamination except itself and the empty set. It is very easy to see that a lamination is minimal if and only if every leaf is dense in  the lamination.  We say that a lamination $\mu$ is an exceptional minimal lamination, if $\mu$ is minimal and does not have interior in $M$.  Thus, the intersection of a transversal with an exceptional minimal lamination is a Cantor set. The following theorem is one of the fundamental results on measured laminations/foliations, see \cite{CC1} for measured foliations.

\begin{theorem}[Theorem 3.2 in Chapter I of \cite{MS}, pp 410]\label{TMS}
Let $\mu$ be a co-dimension one measured lamination in a closed connected 3--manifold $M$, and suppose $\mu\ne M$.  Then, $\mu$ is the disjoint union of a finite number of sub-laminations.  Each of these sub-laminations is of one of the following types:
\begin{enumerate}
\item A family of parallel compact leaves,
\item A twisted family of compact leaves,
\item An exceptional minimal measured lamination.
\end{enumerate} 
\end{theorem}

One can also naturally define the Euler characteristic for measured laminations, see \cite{MS}.  For example, if a measured lamination consists of a family of parallel compact leaves, then its Euler characteristic is equal to the product of the Euler characteristic of a leaf and the total weight.

Using branched surfaces, Morgan and Shalen gave a combinatorial formula for Euler characteristic of measured laminations.  Let $B$ be a branched surface fully carrying a measured lamination $\mu$.  For each branch sector $b$ of the branched surface $B$, one can define a special Euler characteristic $\chi(b)=\chi_{top}(b)-o(b)/4$, where $\chi_{top}(b)$ is the usual Euler characteristic for surfaces and $o(b)$ is the number of corners of $b$, see Definition 3.1 in Chapter II of \cite{MS} pp 424.  Let $w(b)$ be the coordinate (or weight) of $\mu$ at the branched sector $b$.  Then, $\chi(\mu)=\sum w(b)\cdot\chi(b)$ (see Theorem 3.2 in Chapter II of \cite{MS}, pp 424).  The following proposition is easy to prove.

\begin{proposition}\label{PMS}
Let $\mu\subset M$ be a measured lamination with $\chi(\mu)=0$, and let $B$ be a branched surface fully carrying $\mu$.  Suppose $B$ does not carry any 2--sphere.  Then, $B$ fully carries a collection of tori.
\end{proposition}
\begin{proof}
First note that, by Proposition~\ref{Pnonori}, $B$ does not carry any projective plane.  We can add the equation $\sum\chi(b)\cdot w(b)=0$ to the branch equations, and get a new system of linear equations.  By the formula above, every solution to this linear system corresponds to a measured lamination with Euler characteristic 0.  Since all the coefficients are rational numbers and this linear system has a positive solution $\mu$, this linear system must have a positive integer solution.  Hence, $B$ fully carries a collection of closed surfaces with total Euler characteristic 0.  Since $B$ does not carry any closed surface with positive Euler characteristic, each surface in this collection has Euler characteristic 0.  For any Klein bottle $K$ carried by $B$, the boundary of a twisted $I$--bundle over $K$ is a torus carried by $B$.  So, we can get a collection of tori fully carried by $B$.
\end{proof}

The following theorem of Morgan and Shalen is also useful.

\begin{theorem}[Theorem II 5.1 of \cite{MS}, pp 427]\label{TMS2}
Let $B$ be a branched surface that does not carry any surface of positive Euler characteristic.  Let $\mu$ be a measured lamination fully carried by $B$, and suppose every leaf $l$ of $\mu$ has virtually abelian fundamental group.  Then, any measured lamination $\mu'$ carried by $B$ has $\chi(\mu')=0$.
\end{theorem}

An immediate corollary of Theorem~\ref{TMS2} is the following.

\begin{corollary}\label{CMS}
Let $B\subset M$ be a branched surface that does not carry any 2--sphere.  If $B$ fully carries a measured lamination with Euler characteristic 0, then every measured lamination fully carried by $B$ has Euler characteristic 0.
\end{corollary}
\begin{proof}
If $B$ fully carries a measured lamination with Euler characteristic 0, by Proposition~\ref{PMS}, $B$ fully carries a measured lamination consisting of tori.  Hence, Theorem~\ref{TMS2} implies any measured lamination $\mu$ carried by $B$ has $\chi(\mu)=0$.
\end{proof}

Instead of considering the solution space of the system of branch equations, it is more common to consider the projective space, which is usually called the \emph{projective lamination space} (sometimes we also call it the projective solution space).  This is first used by Thurston to study curves and 1-dimensional measured laminations on a surface through the use of train tracks, and it can be trivially generalized to 2-dimensional measured laminations and branched surfaces.  Throughout this paper, we identify the projective lamination space with the set of points $(x_1,\dots,x_N)\in\mathcal{S}(B)$ satisfying $\sum_{i=1}^Nx_i=1$.  We denote the projective lamination space (for the branched surfaces $B$) by $\mathcal{PL}(B)$.  Thus, each rational point of $\mathcal{PL}(B)$ corresponds to a compact surface carried by $B$, and each irrational point corresponds to a measured lamination that contains an exceptional minimal sub-lamination.  By an irrational point, we mean a point in $\mathcal{PL}(B)$ with at least two coordinates are not rationally related.  We may also consider the set of points in $\mathcal{PL}(B)$ corresponding to measured laminations with Euler characteristic 0.  The following proposition is obvious after adding the combinatorial formula of Euler characteristic into the linear system of branch equations.

\begin{proposition}\label{Pcompact}
Let $\mathcal{T}(B)\subset\mathcal{PL}(B)$ be the subset of points corresponding to measured laminations with Euler characteristic 0.  Then $\mathcal{T}(B)$ is a closed and hence compact subset of $\mathcal{PL}(B)$.
\end{proposition}\qed

\section{Measured laminations with Euler characteristic 0}\label{Smin}

The goal of this section is to prove Lemma~\ref{Llocus} which is a certain characterization of measured laminations with Euler characteristic 0.  Lemmas~\ref{Lnodoc} and \ref{Lvan} are also used in \cite{L4}.  The proof involves some basic properties of foliations and laminations, such as the Reeb stability theorem and local stability theorem.  We refer to \cite{CN,CC1,Ta} for more details, see also \cite{GO} for lamination versions of these results.  The Reeb stability theorem basically says that the holonomy along a trivial curve in a leaf must be trivial.  The simplest version of the local stability theorem (for our purpose) basically says that, for any disk $\Delta$ in a leaf, there is a 3--ball neighborhood of $\Delta$ in $M$ whose intersection with the lamination consists of disks parallel to $\Delta$. 

The proof of next lemma is similar to some arguments in section 2 of \cite{L1}.  

\begin{lemma}\label{Lnodoc}
Let $B$ be a branched surface fully carrying a measured lamination $\mu$.  Suppose $\partial_hN(B)$ has no disk component and $N(B)$ does not carry any disk of contact that is disjoint from $\mu$.  Then, $N(B)$ does not carry any disk of contact.
\end{lemma}
\begin{proof}
After some isotopy, we may assume $\partial_hN(B)\subset\mu$ (note that if $\mu$ is a compact surface, we may need to take multiple copies of $\mu$ to achieve this).  For any component $E$ of $\partial_hN(B)$, let $l_E$ be the leaf of $\mu$ containing $E$.  Suppose $l_E-int(E)$ has a disk component $D$.  Note that $D$ is a disk of contact by definition.  Since $\partial_hN(B)$ has no disk component, we may choose $E$ so that $D$ does not contain any component of $\partial_hN(B)$.  Then, after a small isotopy, we can get a disk of contact parallel to $D$ and disjoint from $\mu$.  So, by our hypotheses, $E$ must be an essential sub-surface of $l_E$ and $E$ is not a disk.

After replacing non-orientable leaves by $I$--bundles over these leaves and then deleting the interior of these $I$--bundles (operations 2.1.1--2.1.3 in \cite{G1}), we may assume every leaf of $\mu$ is orientable. After applying these operations to each leaf, we may also assume $\mu$ is nowhere dense \cite{G1}. Suppose there is a disk of contact $D\subset N(B)$.  We may assume $\partial D\subset int(\partial_vN(B))$, $D\cap\mu\subset int(D)$, and $D$ is transverse to each leaf of $\mu$.  Since $\mu$ is a measured lamination, there is no holonomy and every component of $D\cap\mu$ is a circle.  For any circle $\alpha\subset D\cap\mu$, we denote by $\Delta_\alpha$ the disk in $D$ bounded by $\alpha$ and denote the leaf of $\mu$ containing $\alpha$ by $l_\alpha$.  The circle $\alpha$ has two annular collars $A_\alpha^+$ and $A_\alpha^-$ in $l_\alpha$ on the two sides of $\alpha$, where $A_\alpha^+\cap A_\alpha^-=\alpha$ and $A_\alpha^+\cup A_\alpha^-$ is a regular neighborhood of $\alpha$ in $l_\alpha$.  We may assume $A_\alpha^+$, the plus side of $\alpha$, is the one with the property that (after smoothing out the corners) the surface $A_\alpha^+\cup (D-int(\Delta_\alpha))$ is transverse to the $I$--fibers of $N(B)$, hence (after smoothing out the corners) $A_\alpha^-\cup\Delta_\alpha$ is transverse to the $I$--fibers of $N(B)$.  We say that $\alpha$ is of type $I$ if $\alpha$ bounds a disk, denoted by $\Delta_\alpha'$, in $l_\alpha$ and $A_\alpha^+\subset\Delta_\alpha'$, see Figure~\ref{type}(a) for a one-dimension lower schematic picture.  Otherwise, we say $\alpha$ is of type $II$.  Notice that if $\alpha$ is of type $I$, the canonical cutting and pasting of $D$ and $l_\alpha$ at $\alpha$ produce another disk of contact $(D-\Delta_\alpha)\cup\Delta_\alpha'$.   If every circle of $D\cap\mu$ is of type $I$,  we can take the circles of $D\cap\mu$ which are outermost in $D$ and perform the canonical cutting and pasting along these curves.  Then, after some isotopy, we get a disk of contact disjoint from $\mu$.  So, there is at least one type $II$ circle in $D\cap\mu$.

\begin{figure}
\begin{center}
\psfrag{l}{$l_\alpha$}
\psfrag{D}{$D$}
\psfrag{e}{$\Delta_\alpha$}
\psfrag{e'}{$\Delta_\alpha'$}
\psfrag{1}{type $I$}
\psfrag{2}{type $II$}
\psfrag{(a)}{(a)}
\psfrag{(b)}{(b)}
\psfrag{(c)}{(c)}
\includegraphics[width=4in]{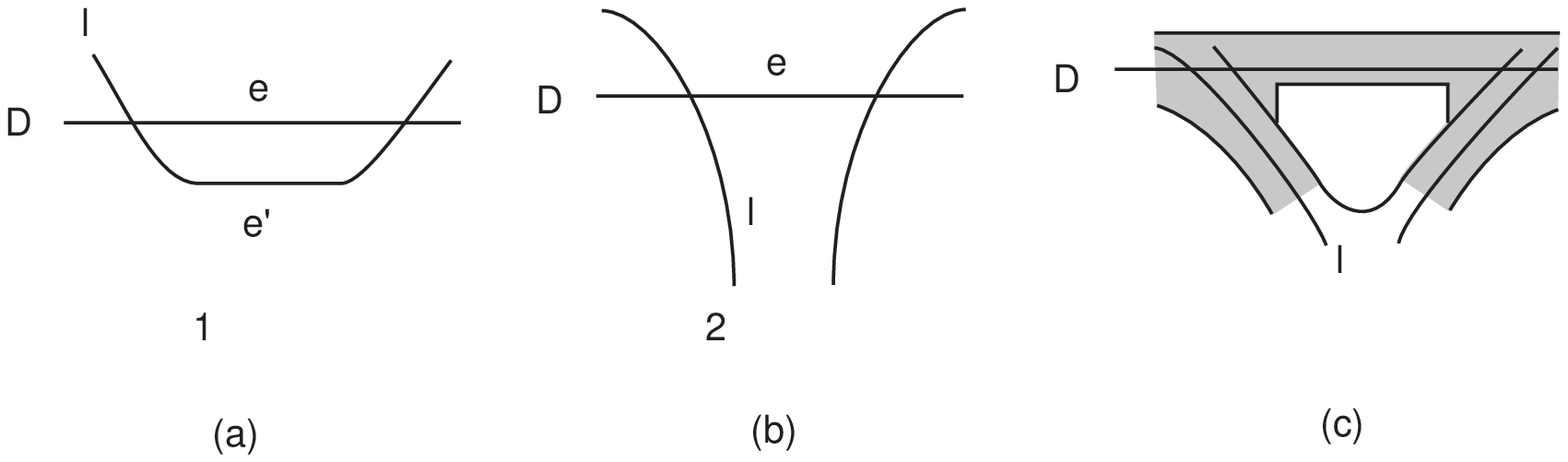}
\caption{}\label{type}
\end{center}
\end{figure}

By the local stability theorem of foliations and laminations, the limit of type $II$ circles of $D\cap\mu$ cannot be a circle of type $I$.  So, we can find a circle $\alpha$ in $D\cap\mu$ such that $\alpha$ is of type $II$ and is innermost in the sense that every circle in $int(\Delta_\alpha)\cap\mu$ is of type $I$.  

Since $\alpha$ is of type $II$, $\alpha$ does not bound a disk in $l_\alpha$ that contains $A_\alpha^+$.  So, one cannot push $\Delta_\alpha$ into $l_\alpha$ along the $I$--fibers of $N(B)$, fixing $\alpha$.  In other words, $\Delta_\alpha$ is not homotopic to a disk in $l_\alpha$ via a homotopy that fixes $\alpha$ and is invariant on each $I$--fiber of $N(B)$.  Therefore, we can find an arc $\beta$ properly embedded in $\Delta_\alpha$, such that one cannot push $\beta$ (fixing $\partial\beta$) into $l_\alpha$ along the $I$--fibers.  Notice that, for any point $x$ near $\partial\beta$, there is a subarc of an $I$--fiber connecting $x$ to a point in $A_\alpha^+\subset l_\alpha$.  We can view $\beta$ as an embedding $\beta:[0,1]\to\Delta_\alpha$.  So, there is a maximal interval $[0,t]$ ($t<1$) such that the arc $\beta([0,t])$ is homotopic to an arc in $l_\alpha$ via a homotopy that fixes $\partial\beta$ and is invariant on each $I$--fiber of $N(B)$.  Thus, for each $\beta(s)$ ($0<s\le t$), there is a subarc $J_s$ of an $I$--fiber such that $\partial J_s$ consists of $\beta(s)$ and a point in $l_\alpha$.  Note that $J_s$ may be degenerate, i.e., $J_s$ may be a single point, in which case $\beta(s)\in l_\alpha\cap\beta$.  We may also regard $J_0$ as the point $\beta(0)$. Since $[0,t]$ is maximal, the arc $J_t$ must contain a vertical arc of $\partial_vN(B)$ (otherwise, one can trivially extend $\beta([0,t])$ along $\beta$ to a longer arc).  This implies that there is an interior point $X$ of $J_t$ such that $X\in\partial_hN(B)\cap\partial_vN(B)$.  We denote the component of $\partial_hN(B)$ containing $X$ by $E_X$ and denote the leaf containing $X$ by $l_X$.  Since $\partial_hN(B)\subset\mu$, $E_X\subset l_X$.  Now, we consider the intersection of $l_X$ and the (singular) triangle $\cup_{s\in[0,t]}J_s$ (the three edges of the triangle are $\beta([0,t])$, $J_t$ and an arc in $l_\alpha$).  As shown in Figure~\ref{type}(c), there must be an arc in $l_X\cap(\cup_{s\in[0,t]}J_s)$ connecting $X$ to a point $\beta(s)$ with $0<s<t$.  Since every circle in $\mu\cap int(\Delta_\alpha)$ is of type $I$, this implies that $X$ lies in a disk of $l_X$ bounded by a type $I$ circle of $\mu\cap int(\Delta_\alpha)$.  Since $X\in\partial E_X$, $E_X$ lies in this disk of $l_X$ bounded by a type $I$ circle of $\mu\cap int(\Delta_\alpha)$, which contradicts our previous conclusion that each component of $\partial_hN(B)$ is a non-disk essential sub-surface of the corresponding leaf.
\end{proof}

%\vspace{5pt}

\noindent\emph{Remark}.  Lemma~\ref{Lnodoc} is true without the hypothesis that $\mu$ is measured.  Suppose $\mu$ is an arbitrary lamination fully carried by $B$.  Then by the Reeb stability theorem, any limiting circle of a spiral in $D\cap\mu$ cannot be of type $I$. So one can proceed as in the proof of Lemma~\ref{Lnodoc} except that a slightly more delicate argument on $D\cap\mu$ is needed in the end.

\begin{definition}\label{Dvan}
Recall that a vanishing cycle (see \cite{CN,GO}) in a foliation $\mathcal{F}$ is an curve $f_0: S^1\to l_0$, where $l_0$ is a leaf in $\mathcal{F}$, and $f_0$ extends to a map $F:[0, 1]\times S^1\to M$ satisfying the following properties.
\begin{enumerate}
\item for any $t\in [0, 1]$, the curve $f_t(S^1)$, defined by $f_t(x)=F(t,x)$, is contained in a leaf $l_t$,
\item for any $x\in S^1$, the curve $t\to F(t,x)$ is transverse to $\mathcal{F}$,
\item $f_0$ is an essential curve in $l_0$, but $f_t$ is null-homotopic in $l_t$ for $t>0$.
\end{enumerate}
We define a slightly different version of vanishing cycle for laminations. Let $\mu$ be a lamination in $M$ and $l_0$ be a leaf.  We call a simple closed curve $f_0:S^1\to l_0$ an \emph{embedded vanishing cycle} in $\mu$ if $f_0$ extends to an embedding $F:[0, 1]\times S^1\to M$ satisfying the following properties.
\begin{enumerate}
\item $F^{-1}(\mu)=C\times S^1$, where $C$ is a closed set of $[0, 1]$, and for any $t\in C$, the curve $f_t(S^1)$, defined by $f_t(x)=F(t,x)$, is contained in a leaf $l_t$, 
\item for any $x\in S^1$, the curve $t\to F(t,x)$ is transverse to $\mu$
\item $f_0$ is an essential curve in $l_0$, but there is a sequence of points $\{t_n\}$ in $C$ such that $\lim_{n\to\infty}t_n=0$ and $f_{t_n}(S^1)$ bounds a disk in $l_{t_n}$ for all $t_n$. 
\end{enumerate}
\end{definition}

\begin{lemma}\label{Lvan}
Let $M$ be a closed orientable and irreducible 3--manifold, and $\mu\subset M$ an exceptional minimal measured lamination.  Suppose $\mu$ is fully carried by a branched surface $B$ and $B$ does not carry any $2$--sphere.  Then, $\mu$ has no embedded vanishing cycle.
\end{lemma}
\begin{proof}
The proof of this lemma is essentially an argument of Novikov.  Novikov showed that \cite{N, CN} if a transversely orientable foliation has a vanishing cycle, then the foliation contains a Reeb component. Note that the $C^2$ assumption in Novikov's original proof is not necessary, see \cite{So} or section 9.3 of \cite{CC2}.  We will use an adaptation of Novikov's argument as in the proof of Lemma 2.8 of \cite{GO}.  

Our proof is based on the proof of Lemma 2.8 of \cite{GO} (pp 54).  So, before we proceed, we briefly describe the argument in \cite{GO}, which shows that a lamination fully carried by an essential branched surface has no vanishing cycle.  In that proof \cite{GO}, the lamination $\lambda$ is fully carried by $N(B)$.  Although the hypothesis of Lemma 2.8 of \cite{GO} is that $B$ is an essential branched surface, the only requirement is that each disk component of $\partial_hN(B)$ is a horizonal boundary component of a $D^2\times I$ region in $M-int(N(B))$.  The first step of the proof in \cite{GO} is to consider $\hat{N}(B)$ which is the union of $N(B)$ and all the $D^2\times I$ regions of $M-int(N(B))$.  So, $\partial_h\hat{N}(B)$ does not contain any disk component. Let $\hat{\mathcal{F}}$ be the associated (singular) foliation of $\hat{N}(B)$ ($\hat{\mathcal{F}}$ is obtained by filling the $I$--bundle regions of $\hat{N}(B)-\lambda$).  The only singularities of $\hat{\mathcal{F}}$ are at $\partial_h\hat{N}(B)\cap\partial_v\hat{N}(B)$.  Then, one simply applies Novikov's argument to the (singular) foliation $\hat{\mathcal{F}}$.  The key of the proof in \cite{GO} is that when one extends the vanishing cycle to a map $F:(0,1]\times D^2\to M$ as in Novikov's argument \cite{N, CN}, the disk $F(\{t\}\times D^2)$ (lying in a leaf of $\hat{\mathcal{F}}$) does not contain any component of $\partial_h\hat{N}(B)$ (since $\partial_h\hat{N}(B)$ has no disk component and $N(B)$ does not carry any disk of contact).  So, the singularities of $\hat{\mathcal{F}}$ never affect Novikov's argument.  Hence, $\hat{\mathcal{F}}$ has a Reeb component.  Note that, by taking a 2-fold cover of $\hat{N}(B)$ if necessary, one can always assume $\hat{\mathcal{F}}$ is transversely orientable.   

Now, we prove Lemma~\ref{Lvan} using the arguments above.  However, since our lamination $\mu$ may be compressible, a disk component of $\partial_hN(B)$ may not correspond to a $D^2\times I$ region of $M-int(N(B))$.  Let $\mathcal{C}$ be the number of components of $M-int(N(B))$ that are not $D^2\times I$ regions.  We assume $\mathcal{C}$ is minimal among all such measured laminations and branched surfaces that satisfy the hypotheses of Lemma~\ref{Lvan} and contain embedded vanishing cycles.

Suppose $\gamma$ is an embedded vanishing cycle in $\mu$. So, $\gamma$ is an essential simple closed curve in a leaf.  There is an embedded vertical annulus $A$ in $N(B)$ containing $\gamma$.  Since $\mu$ is a measured lamination, $\mu$ has no holonomy and we may assume $A\cap\mu$ is a union of parallel circles.  Moreover, by Definition~\ref{Dvan} there is a sequence of circles $\{\gamma_n\}$ in $A\cap\mu$ such that $\lim_{n\to\infty}\gamma_n=\gamma$ and each $\gamma_n$ bounds a disk in $l_n$, where $l_n$ is the leaf of $\mu$ containing $\gamma_n$.  Let $D_n$ be the disk bounded by $\gamma_n$ in $l_n$ and suppose $n$ is sufficiently large.  First note that these $D_n$'s are all on the same side of $A$.  More precisely, for any $D_m$ and $D_n$ (with both $m$ and $n$ sufficiently large), there is a map $\phi: D^2\times I\to M$ such that $\phi(\partial D^2\times I)\subset A$ and $\phi(D^2\times\partial I)=D_m\cup D_n$.  This is because if $D_m$ and $D_n$ are on different sides of $A$, then the union of $D_m\cup D_n$ and the sub-annulus of $A$ bounded by $\gamma_m\cup\gamma_n$ is a 2--sphere $S$, and after a small perturbation, $S$ becomes an immersed 2--sphere carried by $N(B)$.  The canonical cutting and pasting on $S$ can produce an embedded 2--sphere carried by $N(B)$, which contradicts our hypothesis.

Let $A'$ be the sub-annulus of $A$ between $\gamma_m$ and $\gamma_n$. By assuming $m$ and $n$ to be sufficiently large, we may assume $\gamma_m$ and $\gamma_n$ are close to $\gamma$ and hence we may regard $A'$ as a vertical annulus in $N(B)$.  We will show next that every circle in $A'\cap\mu$ bounds a disk in the leaf that contains this circle.  We first consider the generic case: $D_m\cap A'=\partial D_m=\gamma_m$ and $D_n\cap A'=\partial D_n=\gamma_n$.  Since $M$ is irreducible and since $m$ and $n$ are sufficiently large, $D_m\cup A'\cup D_n$ must be an embedded 2--sphere bounding a 3--ball $E=D^2\times I$, where $D^2\times\partial I=D_m\cup D_n$ and $\partial D^2\times I=A'$.  If $E-int(N(B))$ consists of $D^2\times I$ regions, then $\mu\cap E$ is a union of parallel disks with boundary in $A'$.  Conversely, if $\mu\cap E$ consists of parallel disks with boundary in $A'$, then after some splitting, $E-int(N(B))$ becomes a union of $D^2\times I$ regions.  Recall that we have assumed that $\mathcal{C}$, the number of non-$D^2\times I$ regions of $M-int(N(B))$, is minimal among all such measured laminations.  We claim that $E-int(N(B))$ consists of $D^2\times I$ regions.  Otherwise, $\mu\cap E$ must contain non-disk leaves.

By the local stability theorem, the union of all non-disk leaves of $\mu\cap E$ form a sub-lamination of $\mu\cap E$, and we denote this sub-lamination of $\mu\cap E$ by $\lambda$.  So, we can obtain a new measured lamination $\mu'$ by cutting off $\lambda$ from $\mu$ and then gluing back disks along the boundary circles of $\lambda$.  The disks that we glue back are parallel to the disk components of $\mu\cap E$, so we may assume the new lamination $\mu'$ is carried (not fully carried) by $N(B)$.  Moreover, $\mu'$ has a transverse measure induced from that of $\mu$.  Next, we show that $\gamma$ is still an embedded vanishing cycle for $\mu'$.  Let $l_\gamma'$ be the leaf of $\mu'$ containing $\gamma$. By the construction, we only need to show that $\gamma$ is an essential curve in $l_\gamma'$.  Suppose $\gamma$ bounds a disk $\Delta_\gamma$ in $l_\gamma'$.  As $\gamma$ is essential in $\mu$ but trivial in $\mu'$, $\Delta_\gamma\cap E\ne\emptyset$.  Since $\lim_{k\to\infty}\gamma_k=\gamma$ and $\lambda$ is a sub-lamination of $\mu\cap E$, if $k$ is sufficiently large, $D_k\cap E\ne\emptyset$ and  $D_k\cap\lambda\ne\emptyset$, where $D_k$ is the disk bounded by $\gamma_k$ in $\mu$ as above.  This implies that $\overline{D_k-\lambda}$ has a disk component $\Delta$ with $\partial\Delta\subset A'$.  Moreover, after a slight perturbation, the union of $\Delta\cup D_n$ and the sub-annulus of $A'$ bounded by $\partial\Delta\cup\partial D_n$ form an immersed 2--sphere transverse to the $I$--fibers.  After some cutting and pasting, one can obtain an embedded 2--sphere carried by $B$, contradicting our hypothesis.   Therefore, $\gamma$ is still an embedded vanishing cycle for the new measured lamination $\mu'$.  After splitting $N(B)$ along $\mu'\cap E$ and taking sub-branched surfaces, each component of $E-int(N(B))$ becomes a $D^2\times I$ region.  This contradicts our assumption that $\mathcal{C}$ is minimal for $\mu$.  So, in this generic case, every circle of $A'\cap\mu$ bounds a disk in the leaf. 

The non-generic case is very similar.  If $D_m\cap int(A')=\emptyset$ but $D_n\subset int(D_m)$, then we have a map $\phi: D^2\times I\to M$ such that $\phi(\partial D^2\times I)=A'$, $\phi(D^2\times\{0\})=D_m$, $\phi(D^2\times\{1\})=D_n$, and $\phi$ restricted to $D^2\times(0,1)$ is an embedding (this is a standard picture in Novikov's argument on Reeb components, see pp 133 of \cite{CN} for a picture).  So, we can apply the argument above to the (half open) 3--ball $\phi(D^2\times(0,1))$ and the proof is the same. 
If $D_m\cap int(A')\ne\emptyset$, then we can replace $D_n$ and $A'$ by a sub-disk of $D_m$ and a sub-annulus of $A'$ respectively and return to the case that $D_m\cap int(A')=\emptyset$.  Thus, after choosing a sub-annulus of $A$, we may assume that $\gamma\subset\partial A$ and every circle in $\mu\cap(A-\gamma)$ bounds a disk in the corresponding leaf.

Let $D_m$, $D_n$, $A'$ and $E=D^2\times I$ be as above.  By the arguments above, $\mu\cap E$ consists of parallel disks.  As $A'\subset int(N(B))$, if the 3--ball $E$ contains some components of $M-N(B)$, then we can split $N(B)$ in $E$ so that $E-int(N(B))$ consists of $D^2\times I$ regions.  We can perform this splitting in all possible such 3--balls $E$, and this is a finite process since $|\partial_hN(B)|$ is bounded.  Let $\hat{N}(B)$ be the union of $N(B)$ (after the splitting in the 3--balls $E$ above) and all the $D^2\times I$ regions of $M-int(N(B))$.  The remaining proof is the same as the proof of Lemma 2.8 of \cite{GO}.  We can extend $\mu$ to a (singular) foliation $\hat{\mathcal{F}}$ in $\hat{N}(B)$.  By our construction above, for any disk $\Delta$ (in a leaf) bounded by a circle of $\mu\cap A$, $\Delta$ does not contain any component of $\partial_h\hat{N}(B)$, and hence $\Delta$ does not meet the singularity of $\hat{\mathcal{F}}$.  So, we can apply Novikov's argument to $\hat{\mathcal{F}}$ as in the proof of Lemma 2.8 of \cite{GO}, and conclude that $\hat{\mathcal{F}}$ contains a Reeb component and hence has non-trivial holonomy. Since $\hat{\mathcal{F}}$ is obtained by filling the $I$--bundle regions of $\hat{N}(B)-\mu$, this implies that $\mu$ has non-trivial holonomy and is not a measured lamination.
\end{proof}
\begin{remark}
One can apply Novikov's argument directly to laminations without using the (singular) foliation $\hat{\mathcal{F}}$.  Moreover, the assumption that $B$ does not carry any 2--sphere seems unnecessary.  One can prove Lemma~\ref{Lvan} (without the 2--sphere assumption) using the argument of Imanishi, which says that any 2-dimensional phenomenon like the Reeb foliation implies that the foliation/lamination has non-trivial holonomy, see \cite{L3} for an interpretation of Imanishi's argument using branched surfaces.  
\end{remark}

\begin{lemma}\label{Llocus}
Let $M$ be a closed orientable and irreducible 3--manifold and suppose $M$ is not $T^3=S^1\times S^1\times S^1$.  Let $\mu\subset M$ be an exceptional minimal measured lamination with Euler characteristic 0.  Suppose $\mu$ is fully carried by a branched surface $B$, and $B$ does not carry any 2--sphere.  Then, there is a branched surface $B'$, obtained by splitting $B$ and taking sub-branched surfaces, such that $B'$ fully carries $\mu$, the branch locus $L'$ of $B'$ has no double point, and $B'-L'$ consists of annuli and M\"{o}bius bands.  
\end{lemma}
\begin{proof}
Suppose every leaf of $\mu$ is a plane.  After eliminating all the disks of contact of $N(B)$ that are disjoint from $\mu$, we have that $\partial_hN(B)$ consists of disks.  So there is no monogon and $\mu$ is an essential lamination.  By a Theorem in \cite{G7} (see also Proposition 4.2 of \cite{L1}), $M\cong T^3$.  

So, at least one leaf of $\mu$ is not a plane.  Let $\gamma$ be an essential simple closed curve in a non-plane leaf of $\mu$.  Then, there is a vertical annulus $A$ in $N(B)$ containing $\gamma$.  Since $\mu$ is a measured lamination and so has no holonomy, we may assume $A\cap\mu$ consists of circles parallel to $\gamma$.  If $\mu$ contains a plane leaf $L$, since every leaf of $\mu$ is dense, $L\cap A$ contains an infinite sequence of circles whose limit is $\gamma$.  Each circle in $L\cap A$ bounds a disk in the plane $L$, so $\gamma$ is an embedded vanishing cycle and we get a contradiction to Lemma~\ref{Lvan}.  Thus, $\mu$ contains no plane leaf at all.

After some isotopy, we may assume $\partial_hN(B)\subset\mu$.  Since $\mu$ contains no plane leaf, for every component $S$ of $\partial_hN(B)$, we can split $N(B)$ along $\mu$ so that $S$ contains an essential curve in the leaf that contains $S$.  If there is a disk of contact in $N(B)$ disjoint from $\mu$, then we can trivially eliminate the disk of contact by splitting $B$.   After these splittings, each component $S$ of $\partial_hN(B)$ becomes an essential non-disk sub-surface of the leaf that contains $S$.  By Lemma~\ref{Lnodoc}, $N(B)$ does not contain any disk of contact.

Next, we show that each component of $\partial_hN(B)$ must be an annulus.  By Proposition~\ref{PMS}, $B$ fully carries a collection of tori $T$.  After some isotopy and taking multiple copies of $T$, we may assume $\partial_hN(B)\subset T$.  If a component $S$ of $\partial_hN(B)$ is not an annulus, since no component of $\partial_hN(B)$ is a disk, there must be a boundary component of $S$ bounding a disk $D$ in $T-int(S)$, and $D$ is a disk of contact by definition.  Since we have assumed that $\partial_hN(B)\subset\mu$ after isotopy, if a component of $\partial_hN(B)$ is a closed surface, $\mu$ must contain a closed surface, a contradiction to the hypothesis that $\mu$ is exceptional minimal. Thus, $\partial_hN(B)$ does not contain a torus and $\partial_hN(B)$ must consist of annuli.  

If a leaf of $\mu$ has non-zero Euler characteristic, then we can split $N(B)$ by ``blowing air" into $N(B)-\mu$ so that a component of $\partial_hN(B)$ is an essential sub-surface of a leaf and has negative Euler characteristic. So, the argument above implies that each leaf of $\mu$ is either an infinite annulus or an infinite M\"{o}bius band.  

Let $\eta$ be an essential simple closed curve in $int(\partial_hN(B))$ and let $A_\eta$ be a vertical annulus in $N(B)$ containing $\gamma$.  So, we may assume $A_\eta\cap\mu$ is a union of parallel circles.  By Lemma~\ref{Lvan}, $\eta$ is not an embedded vanishing cycle, hence we can choose the vertical annulus $A_\eta$ so thin that every circle of $A_\eta\cap\mu$ is an essential curve in a leaf of $\mu$.  Since every leaf is dense in $\mu$, each leaf must intersect $A_\eta$.  Moreover, the limit of each end of any leaf  is a sub-lamination of $\mu$ and hence is the whole of $\mu$. So, each end of any leaf of $\mu$ must intersect $A_\eta$.  

After some splittings, we may also assume $|M-N(B)|$ is minimal among all such branched surfaces.  Let $S$ be a component of $\partial_hN(B)$ and $l_S$ be the leaf of $\mu$ containing $S$.  We first point out that $l_S$ must be an orientable surface.  To see this, for any point $x$ in any leaf $l$ and for any transversal $\delta_x$ containing $x$, since every leaf is dense, $x$ is always an accumulation point in $\delta_x\cap\mu$.  Since $M$ is orientable, if $l$ is a non-orientable surface, $x$ must be a limit point (of $\delta_x\cap\mu$) in both components of $\overline{\delta_x-x}$.  However, if $x\in\partial_hN(B)$, $x$ can only be a limit point on one side of $\delta_x$.  So, $l_S$ must be orientable and hence $l_S$ is an infinite annulus.  Since both ends of $l_S$ intersect $A_\eta$ and no circle in $A_\eta\cap l_S$ bounds a disk in $l_S$, there is an annulus in each component of $l_S-int(S)$ connecting $\partial S$ to $A_\eta$.  Therefore, we can find an annulus $A_S$ in $N(B)-\mu$, transverse to the $I$--fibers and with one boundary circle in $A_\eta$ and the other boundary circle in a component of $\partial_vN(B)$. Moreover, $A_S$ is parallel to a sub-annulus  of $l_S-int(S)$ above.  We can split $N(B)$ by deleting a fibered neighborhood of $A_S$ from $N(B)$. Note that since we have assumed $|M-N(B)|$ is minimal, the branched surface after this splitting still fully carries $\mu$ and satisfies all the previous properties.  Since both components of $l_S-int(S)$ contain such annuli, we can find such an annulus in $N(B)$ connecting $A_\eta$ to each component of $\partial_vN(B)$.   By deleting a small neighborhood of these annuli from $N(B)$, we can split $N(B)$ into $N(B')$ which is a fibered neighborhood of another branched surface $B'$ and $N(B')$ satisfies all the previous properties.  Since the splittings are along the annuli connecting $\partial_vN(B)$ to $A_\eta$, each component of $\partial_vN(B')$ lies in a small neighborhood of $A_\eta$ and is parallel to a sub-annulus of $A_\eta$.  Thus, after a small perturbation in a neighborhood of $A_\eta$, we may assume $\pi(\partial_vN(B'))=L'$ is a collection of disjoint circles in $B'$, where $\pi:N(B')\to B'$ is the map collapsing each $I$--fiber to a point.  So, the branched surface $B'$ satisfies all the requirements in Lemma~\ref{Llocus}.
\end{proof}

\section{Normal tori and 0--efficient triangulations}\label{S0eff}

Let $F$ be an embedded surface in $M$ and suppose $M$ has a triangulation $\mathcal{T}$.  We use $\mathcal{T}^{(i)}$ to denote the $i$--skeleton of $\mathcal{T}$.  After some isotopy, we may assume $F$ does not contain any vertices of the triangulation and $F$ is transverse to $\mathcal{T}^{(1)}$ and $\mathcal{T}^{(2)}$. If $F$ is not a normal surface, we can try to normalize $F$ using the following two types of normal moves.  After these normal moves, $F$ consists of normal surfaces and possibly some trivial 2--spheres in 3--simplices.  Note that if $F$ is incompressible, then the two normal moves are isotopies and there are no such 2--spheres. We refer to section 3.1 of \cite{JR} for more detailed descriptions.

\begin{operation}\label{o1}
Suppose $F$ is compressible in a 3--simplex, then there are two cases. The first case is that, for a 2--simplex $\Delta$, $F\cap\Delta$ contains circles.  Let $c$ be a circle of $F\cap\Delta$ innermost in $\Delta$.  If $c$ is a trivial circle in $F$, then the two disks bounded by $c$ in $F$ and $\Delta$ form a 2--sphere bounding a 3--ball.  So, we can perform an isotopy on $F$ pushing the disk across this 3--ball and reduce the number of circles in $F\cap\Delta$.  If $c$ is non-trivial in $F$, the disk bounded by $c$ in $\Delta$ is a compressing disk for $F$ and we can compress $F$ along this compressing disk.  The latter operation increases the Euler characteristic of $F$ by $2$.  The second case is that $F\cap\Delta$ contains no circle but $F$ is compressible in the interior of a 3--simplex.  Similar to the first case, we can either compress $F$ in the interior of the 3--simplex increasing the Euler characteristic, or perform some isotopy reducing the intersection of $F$ with the 2--skeleton.
\end{operation}

\begin{operation}\label{o2}
This operation is an isotopy on $F$.  For any 3--simplex $X$, if $F\cap X$ is incompressible in $X$ and a component of $F\cap X$ intersects an edge of $X$ in more than one point, then one can find a $\partial$--compressing disk $D\subset X$ with $\partial D$ consisting of an arc in $F$ and an arc in an edge (technically $D$ is a $\partial$-parallel disk in $X$).  We can perform an isotopy by pushing $F$ along $D$ across this edge.  This operation reduces the weight of $F$ by two.
\end{operation}

In this section, we will assume the triangulation $\mathcal{T}$ is a $0$--efficient triangulation.  A triangulation of $M$ is said to be $0$--efficient if the triangulation has only one vertex and the only normal 2--sphere in $M$ is the boundary sphere of a closed neighborhood of this vertex.  In \cite{JR}, Jaco and Rubinstein showed that, if $M$ is irreducible and not a lens space, then $M$ admits a $0$--efficient triangulation.  In fact, given any triangulation of $M$, there is an algorithm to collapse this triangulation into a 0--efficient one.  One of the most useful techniques in \cite{JR} is the so-called \emph{barrier surfaces} or \emph{barriers}.  We will briefly explain a special case of barriers used in our proof, see section 3.2 of \cite{JR} for more details.

Let $F$ be a compact embedded normal surface in $M$.  If we cut $M$ open along $F$, we get a manifold with boundary, denoted by $\overline{M-F}$, with an induced cell decomposition.  Let $S$ be a properly embedded normal surface in $\overline{M-F}$ with respect to the induced cell decomposition.  $F\cup S$ is a 2--complex in $M$.  Now we consider the surface $\partial\eta(F\cup S)$ in $M$, where $\eta(F\cup S)$ is the closure of a small neighborhood of $F\cup S$.  The surface $\partial\eta(F\cup S)$ may not be normal and we can use the operations~\ref{o1} and \ref{o2} to normalize $\partial\eta(F\cup S)$.  Then, by \cite{JR}, $F\cup S$ forms a ``barrier" for these normalizing operations.  More precisely, one can perform operations~\ref{o1} and \ref{o2} on $\partial\eta(F\cup S)$ totally in the 3--manifold $M-int(\eta(F\cup S))$ and get a normal surface (with respect to the triangulation of $M$) plus possible trivial 2--spheres in some tetrahedra.  Note that it is possible that, after these operations, $\partial\eta(F\cup S)$ vanishes, i.e. becomes a collection of trivial 2--spheres in some tetrahedra. 

Since every normal 2--sphere in a 0--efficient triangulation is vertex-linking, it is easy to use the barrier technique to derive some nice properties of normal tori with respect to a 0--efficient triangulation.  Lemmas~\ref{Ltorus}, \ref{Lannulus} and Corollary~\ref{Cklein} are well-known to people who are familiar with 0--efficient triangulations.

\begin{lemma}\label{Ltorus}
Suppose $M$ is irreducible and atoroidal and $M$ is not a lens space.  Let $T$ be a normal torus with respect to a 0--efficient triangulation of $M$.  Then, we have the following.
\begin{enumerate}
\item $T$ bounds a solid torus in $M$.  
\item Let $N$ be the solid torus bounded by $T$. Then, $M-int(N)$ is irreducible and $T$ is incompressible in $M-int(N)$.
\end{enumerate}
\end{lemma}
\begin{proof}
As $M$ is irreducible and atoroidal, $T$ is compressible and separating.  Let $D$ be a compressing disk for $T$.  Then, we can choose $D$ so that $D$ is normal with respect to the induced cell decomposition of $\overline{M-T}$.  Hence, $T\cup D$ forms a barrier.  Note that $\partial\eta(T\cup D)$ has a 2--sphere component $S$ and $S$ bounds a 3--ball $E_S$ in $M$.  If this 3--ball $E_S$ lies in the complement of $T\cup D$, then $T$ bounds a solid torus, otherwise $E_S$ contains $T$ and $T$ bounds a ball with a knotted hole.

Since $T\cup D$ forms a barrier, we can perform Operations \ref{o1} and \ref{o2} to normalize $S$ in the complement of $T\cup D$.  Note that Operation~\ref{o2} is an isotopy.  If Operation~\ref{o1} occurs,  since $S$ is a 2--sphere, Operation~\ref{o1} on $S$ is also an isotopy.  Therefore, we can isotope $S$ in $M-T\cup D$ either to a normal 2--sphere or to a 2--sphere in a 3--simplex.  Since the only normal 2--sphere is the vertex-linking one and the normal torus $T$ cannot lie in a small neighborhood of the vertex, $T$ must lie outside the 3--ball bounded by $S$. Hence, $T$ must bound a solid torus.

If $T$ is compressible in the complement of this solid torus $N$, then we have a compressing disk outside the solid torus.  We can use the union of $T$ and this compressing disk as a barrier and the argument above implies that $T$ bounds a solid torus on the other side, which means $M$ is a lens space and contradicts our hypotheses.  If $M-int(N)$ is reducible, then there is an essential normal 2--sphere in $M-int(N)$. Since the only normal 2--sphere is the vertex-linking one and bounds a 3--ball, we also get a contradiction as before. 
\end{proof}

\begin{corollary}\label{Cklein}
Suppose $M$ is a closed, orientable, irreducible and atoroidal 3--manifold and $M$ is not a small Seifert fiber space.  Then, $M$ does not contain any normal projective plane or normal Klein bottle with respect a 0--efficient triangulation.
\end{corollary}
\begin{proof}
If $M$ contains a normal projective plane $P$, then a closed neighborhood of $P$ in $M$, $\eta(P)$, is a twisted $I$--bundle over $P$, and $\partial\eta(P)$ is a normal 2--sphere.  Since the only normal 2--sphere in $M$ is the vertex-linking one, this implies $M$ is $\mathbb{R}P^3$.

If $M$ contains a normal Klein bottle $K$, then $\eta(K)$ is a twisted $I$--bundle over $K$ and $\partial\eta(K)$ is a normal torus.  Since every normal torus bounds a solid torus in $M$, $M$ is the union of a solid torus and a twisted $I$--bundle over a Klein bottle, which implies that $M$ is a Seifert fiber space.
\end{proof}

\begin{lemma}\label{Lannulus}
Suppose $M$ is closed, orientable, irreducible and atoroidal and suppose $M$ is not a small Seifert fiber space.  Let $T$ be a normal torus with respect to a 0--efficient triangulation of $M$, and let $N$ be the solid torus bounded by $T$.  Suppose $A$ is an annulus properly embedded in $M-int(N)$ and $\partial A$ is a pair of essential curves in $T$.  Suppose $A$ is normal with respect to the induced cell decomposition of $M-int(N)$.  Then, the following are true.
\begin{enumerate}
\item each component of $\partial\eta(N\cup A)$ bounds a solid torus in $M$,
\item one component of $\partial\eta(N\cup A)$ bounds a sold torus in $M-int(\eta(N\cup A))$ and the other component of $\partial\eta(N\cup A)$ bounds a solid torus containing $N\cup A$.
\item If $\partial A$ is a pair of meridian curves for the solid torus $N$, then $A$ is $\partial$--parallel in $M-int(N)$.
\end{enumerate}
\end{lemma}
\begin{proof}
Since $\partial A$ is essential in $T$, $\partial\eta(N\cup A)$ consists of two tori in $M-N$.  Let $T_1$ be a component of $\partial\eta(N\cup A)$.  The torus $T_1$ may not be normal, but $T\cup A$ forms a barrier and we can perform Operations \ref{o1} and \ref{o2} to normalize $T_1$ in $M-N\cup A$.  During the normalization process, every step is an isotopy unless in Operation~\ref{o1}, there is a circle in $T_1\cap\Delta$ ($\Delta$ is a 3--simplex) bounding a compressing disk $D$ in $\Delta$.  If this happens, we compress $T_1$ along $D$ as in Operation~\ref{o1} and change $T_1$ into a 2--sphere $T_1'$.  After the compression, similar to the proof of Lemma~\ref{Ltorus}, we can isotope the 2--sphere $T_1'$ either to a normal 2--sphere or into a 3--simplex.  As in the proof of Lemma~\ref{Ltorus}, $T_1'$ must bound a 3--ball in $M-N\cup A$.  Since $N$ and the compressing disk $D$ are on different sides of $T_1$, similar to the proof of Lemma~\ref{Ltorus}, $T_1$ must bound a solid torus in $M-N\cup A$.  If the compression operation never happens, then we can isotope the torus $T_1$ either to a normal torus, in which case $T_1$ bounds a solid torus by Lemma~\ref{Ltorus}, or into a 3--simplex.  If $T_1$ can be isotoped into a 3--simplex, then we have a 3--ball containing $T_1$ and disjoint from $N\cup A$.  This is impossible because the region between $T_1$ and $N\cup A$ is a product.  Thus, each torus in $\partial\eta(N\cup A)$ must bound a solid torus in $M$. 

Let $T_1$ and $T_2$ be the two tori in $\partial\eta(N\cup A)$, and let $E_1$ and $E_2$ be the two components of $M-int(\eta(N\cup A))$ bounded by $T_1$ and $T_2$ respectively.  So, $\partial E_i=T_i$ and each $T_i$ bounds a solid torus in $M$.  If both $E_1$ and $E_2$ are solid tori, then $M$ is a union of $T\cup A$ and 3 solid tori, which implies that either $M$ is a small Seifert fiber space or $M$ is reducible.  Thus, at lease one $E_i$ is not a solid torus.  Suppose $E_1$ is not a solid torus.  Since $T_1$ bounds a solid torus in $M$, $M-int(E_1)$ is a solid torus containing $N\cup A$.  Moreover, by Lemma~\ref{Ltorus}, $M-N$ is irreducible and hence $E_1$ is irreducible.  Since $E_1$ is not a solid torus and $E_1$ is irreducible, $T_1$ must be incompressible in $E_1$.  We claim that $E_2$ must be a solid torus.  Suppose $E_2$ is not a solid torus either. Then the argument above implies that $T_2$ is incompressible in $E_2$.  Let $D_i$ be a meridian disk of the solid torus $M-int(E_i)$ ($i=1,2$).  We first show that at least one of $D_1$ and $D_2$ is properly embedded in $M-int(E_1\cup E_2)$ after isotopy. Suppose $D_1\cap E_2\ne\emptyset$.  Since $M$ and $E_2$ are irreducible, an isotopy can eliminate curves in $D_1\cap T_2$ that are trivial in $T_2$ and innermost in $D_1$.  Since $T_2$ is incompressible in $E_2$, if $D_1\cap T_2\ne\emptyset$ after this isotopy, the subdisk $\Delta$ of $D_1$ bounded by an innermost circle of $D_1\cap T_2$ in $D_1$ is a meridian disk of the solid torus $M-int(E_2)$.  By choosing $D_2$ to be $\Delta$, we have that $D_2$ is properly embedded in $M-int(E_1\cup E_2)$ and clearly $D_2\cap T_1=\emptyset$. Now suppose $D_2$ is properly embedded in $M-int(E_1\cup E_2)$.  Since $M-int(E_2)$ is a solid torus, by compressing $T_2$ along $D_2$, we get a 2--sphere $S_2$ bounding a 3--ball and the 3--ball contains $E_1$.  As $E_1$ lies in this 3--ball, this means that the 2--sphere $S_2$ lies in the solid torus $M-int(E_1)$ and hence bounds a 3--ball in the solid torus $M-int(E_1)$.  Hence $M$ must be $S^3$, a contradiction.  So exactly one of $E_1$ and $E_2$ is a solid torus and part 2 of the lemma holds.

Suppose $\partial A$ is a pair of meridian curves for $N$.  By part 2 of the lemma, $\partial A$ bounds an annulus $A'\subset T$ such that $A\cup A'$ bounds a solid torus $N'$ in $M-int(N)$.  Moreover, each circle in $\partial A$ bounds a meridian disk of $N$.  Since $M$ is not a lens space and $M$ is irreducible, $\partial A$ must be longitudes for the solid torus $N'$.  Thus, $A$ is isotopic to $A'$ (fixing $\partial A$) in $N'$, and part 3 holds.
\end{proof}

Let $B$ be a branched surface in $M$ constructed by gluing normal disks together near the 2--skeleton, as in \cite{FO} and section~\ref{Spre}. By this construction, every surface carried by $B$ is a normal surface.   Let $T$ be a normal surface fully carried by $B$, and we suppose $T\subset N(B)$ and $\partial_hN(B)\subset T$.  So, $\partial_vN(B)$ is a union of annuli properly embedded in $\overline{M-T}$.  By the construction of $B$, after a small perturbation and eliminating disks of contact, we may assume $\partial_vN(B)$ is normal with respect to the induced cell decomposition of $\overline{M-T}$.  

\begin{lemma}\label{L0eff}
Let $M$ be a closed orientable irreducible and atoroidal 3--manifold with a 0--efficient triangulation. Suppose $M$ is not a Seifert fiber space.  Let $B$ be a branched surface as above, i.e., $B$ is obtained by gluing together normal disks, $B$ fully carries a normal surface $T$ with $\partial_hN(B)\subset T$, and $\partial_vN(B)$ is normal with respect to the induced cell decomposition of $\overline{M-T}$.   Suppose the branch locus $L$ of $B$ does not have any double point, $B-L$ consists of annuli and M\"{o}bius bands, and every component of $\partial_hN(B)$ is an annulus.   Then,
\begin{enumerate}
\item $\partial_hN(B)$ is incompressible in $M-int(N(B))$,
\item some component of $\partial N(B)$ bounds a solid torus in $M$ that contains $N(B)$,
\item $M-int(N(B))$ contains a $monogon\times S^1$ region.
\end{enumerate}
\end{lemma}
\begin{proof}
By the hypotheses, any closed surface carried by $B$ is a normal surface with Euler characteristic 0.  By Corollary~\ref{Cklein}, $M$ does not contain any normal Klein bottle.  So, every closed surface carried by $B$ consists of normal tori.  Let $T=\cup_{i=1}^mT_i$ be a collection of disjoint normal tori fully carried by $N(B)$, where each $T_i$ is a component of $T$, and we may assume $\partial_hN(B)\subset T$.  Hence $\partial_vN(B)$ is a collection of annuli properly embedded in $\overline{M-T}$, whose boundary consists of essential curves in $T$.  Moreover, $\partial_vN(B)$ is normal with respect to the induced cell decomposition of $\overline{M-T}$.

By the hypotheses, every component of $\partial N(B)$ is a torus.  Similar to the proof of Lemma~\ref{Lannulus}, $T\cup\partial_vN(B)$ form a barrier, and each component of $\partial N(B)$ bounds a solid torus in $M$. Let $E_1,\dots, E_n$ be the components of $M-int(N(B))$.  Each $\partial E_i$ bounds a solid torus in $M$.  

Suppose $E_i$ is not a solid torus, then $M-int(E_i)$ is a solid torus that contains $N(B)$ and $T$.  Moreover, by the proof of Lemma~\ref{Lannulus}, $E_i$ is irreducible.  Since $E_i$ is not a solid torus, this implies that $\partial E_i$ is incompressible in $E_i$.  Thus, similar to the proof of Lemma~\ref{Lannulus}, for any two components $E_i$ and $E_j$, at least one must be a solid torus.  This implies that at most one component of $M-int(N(B))$ is not a solid torus.  

Next, we show that $\partial_hN(B)$ is incompressible in $M-int(N(B))$.  The basic idea of the proof is that, if $\partial_hN(B)$ is compressible in $M-int(N(B))$, one can construct a solid torus bounded by a new normal torus carried by $N(B)$, and one can use the compressing disk of $M-int(N(B))$ to obtain a compressing disk of this new normal torus outside this solid torus, which contradicts part 2 of Lemma~\ref{Ltorus}.  This solid torus is constructed by joining two $monogon\times S^1$ regions of $M-int(N(B))$.

Let $N_i$ be the solid torus bounded by $T_i$. Suppose $\partial_hN(B)$ is compressible in $M-int(N(B))$ and let $D$ be a compressing disk.  We may suppose $\partial D\subset\partial_hN(B)$ lies in $T_1$ and by Lemma~\ref{Ltorus}, $D$ is a meridian disk of the solid torus $N_1$ bounded by $T_1$.  

Let $H=\overline{N(B)-T}$. Since $\partial_hN(B)\subset T$, $H$ is a collection of $annuli\times I$ and twisted $I$--bundles over M\"{o}bius bands.  $\partial H$ consists of two parts, the horizontal boundary $\partial H\cap T$ and the vertical boundary $\partial H\cap\partial_vN(B)$.  We denote the horizontal boundary of $H$ by $\partial_hH$ and the vertical boundary of $H$ by $\partial_vH$ ($\partial_vH=\partial_vN(B)$).  By the hypotheses, $\partial_hH$ consists of essential annuli in $T$.   Since no component of $\partial_hN(B)$ is a torus, both $N_1$ and $M-int(N_1)$ contain some components of $H$.  

So, there must be a component of $\partial_vN(B)$, say $V$, properly embedded in $N_1$.  By our assumptions, if $V$ is not $\partial$--parallel in $N_1$, then $V$ can be obtained by attaching a knotted tube to a pair of compressing disks of $N_1$.  This implies that a component of $\overline{N_1-V}$, say $\Sigma$, is a 3--ball with a knotted hole.  So, $\partial\Sigma$ is a torus incompressible in $\Sigma$.  Since $T_1\cup V$ forms a barrier, we can use Operations~\ref{o1} and \ref{o2} to isotope $\partial\Sigma$ into a normal torus in $\Sigma$.  However, by Lemma~\ref{Ltorus}, $M-int(\Sigma)$ must be a solid torus.  Since $\Sigma$ is a ball with a knotted hole, this implies that $M$ is $S^3$.  Therefore, each component of $\partial_vN(B)$ in $N_1$ must be $\partial$--parallel in $N_1$.  This implies that there must be a $monogon\times S^1$ region of $M-int(N(B))$ in $N_1$.  We denote this $monogon\times S^1$ region by $J_1$.  So, $\partial J_1$ consists of an annulus in $T_1$ and a component of $\partial_vN(B)$, and $J_1\cap D=\emptyset$ ($D$ is the compressing disk above).

Now, we consider the components of $H$ that lie in $M-int(N_1)$. The simplest case is that there is a component of $H$, say $H_1$, in $M-int(N_1)$ with its horizontal boundary totally in $T_1$. By the construction, the vertical boundary of $H_1$ consists of annuli properly embedded in $M-int(N_1)$.  Since the branch locus $L$ has no double point and $\partial_hN(B)$ is compressible in $N_1$, the boundary curves of $\partial_vH_1$ are meridian curves in $\partial N_1$.  By part 3 of Lemma~\ref{Lannulus}, each annulus in $\partial_vH_1$ is $\partial$--parallel in $M-int(N_1)$.  So, there is also a $monogon\times S^1$ region $J_2$ of $M-int(N(B))$ in $M-int(N_1)$ with $\partial J_2$ consisting of an annulus in $T_1$ and a component of $\partial_vN(B)$.  We denote the component of $\partial_vN(B)$ in $\partial J_i$ by $V_i$ ($i=1,2$).  Within a small neighborhood of $T_1$ in $N(B)$, we can find an annulus $A\subset N(B)$ connecting $V_1$ to $V_2$ and transverse to the $I$--fibers of $N(B)$.  The union of $\partial J_1-V_1$, $\partial J_2-V_2$ and two parallel copies of $A$ form a torus $T_J$ carried by $N(B)$, and $T_J$ bounds a solid torus $N_J$ which is the union of $J_1$, $J_2$ and a product neighborhood of $A$. By the hypothesis on $B$, $T_J$ is a normal torus.  However, since the boundary of $\partial_vH_1$ consists of meridian curves of $\partial N_1$, a meridian disk of $N_1$ gives rise to a compressing disk for the torus $T_J$ in $M-int(N_J)$.  This contradicts part 2 of Lemma~\ref{Ltorus}.

Suppose there is a component $H_2$ of $H$ with one horizontal boundary component in $T_1$ and the other horizontal boundary component in $T_2$.  Then, by our assumption on the meridian curves, the union of a vertical annulus of $H_2$ and a meridian disk of $N_1$ form a compressing disk for $T_2$.  By part 2 of Lemma~\ref{Ltorus}, we must have $N_1\subset N_2$.  Suppose $N_1\subset\dots\subset N_k$ are a maximal collection of nested solid tori, such that there is a component of $H$ between each pair of tori $T_i\cup T_{i+1}$, same as the $H_2$ above.  Since $k$ is maximal, there must be a component of $H$ in $M-int(N_k)$ with horizonal boundary totally in $T_k$.  As before, there is a $monogon\times S^1$ region $J_2$ of $M-int(N(B))$ in $M-int(N_k)$ with $\partial J_2$ consisting of an annulus in $T_k$ and a component of $\partial_vN(B)$.  By assembling annuli in the $T_i$'s ($i=1,\dots,k$) and annuli in those components of $H$ between the tori $T_i\cup T_{i+1}$, we can construct an annulus $A\subset N(B)$, such that $A$ connects $J_1$ to $J_2$ as before and $A$ is transverse to the $I$--fibers of $N(B)$.  Similarly, we can form a torus $T_J$ bounding a solid torus $N_J$, and $N_J$ is the union of $J_1$, $J_2$ and a product neighborhood of $A$.  Moreover, a meridian disk of $N_1$ gives rise to a compressing disk for $T_J$ in $M-int(N_J)$, and we get a contradiction to part 2 of Lemma~\ref{Ltorus}.  This proves that $\partial_hN(B)$ is incompressible in $M-int(N(B))$.

By the hypotheses, $N(B)$ is a Seifert fiber space, and the Seifert fibration restricted to each annulus $\partial_hN(B)$ or $\partial_vN(B)$ is the standard foliation by circles.  If every component of $M-int(N(B))$ is a solid torus, then since $\partial_hN(B)$ is incompressible in $M-int(N(B))$, $M$ is a Seifert fiber space.  Therefore, by the conclusion before, exactly one component $E_i$ of $M-int(N(B))$ is not a solid torus, and $M-int(E_i)$ is a solid torus containing $N(B)$.  

Let $N_1$ be  an innermost solid torus.  By the argument before, each component of $\partial_vN(B)\cap N_1$ is $\partial$--parallel in $N_1$.  This implies that there is a $monogon\times S^1$ region in $N_1$, and part 3 of the lemma holds.
\end{proof}

\section{Splitting branched surfaces, the torus case}\label{Storus}

A main technical part of this paper is to show that, if a branched surface $B$ carries a sequence of Heegaard surfaces $\{S_n\}$ and a measured lamination $\mu$ with $\chi(\mu)=0$, then one can split $B$ into a collection of branched surfaces, such that each $S_n$ is carried by a branched surface in this collection and no branched surface in this collection carries $\mu$.  In this section, we consider the case that $\mu$ is a torus, and we prove the case that $\mu$ is an exceptional minimal lamination in the next section.  The goal of this section is to prove Lemma~\ref{Ltorus2}.

Let $B$ be a branched surface carrying a sequence of closed orientable surfaces $\{S_n\}$.  Suppose $\mu$ is a lamination carried (but may not be fully carried) by $B$.  By section~\ref{Spre}, there is a sub-branched surface of $B$, denoted by $B_\mu$, fully carrying $\mu$.  We may consider $N(B_\mu)\subset N(B)$ with compatible $I$--fiber structure.  Let $D\subset N(B_\mu)\subset N(B)$ be a disk transverse to the $I$--fibers.  We call $D$ a \emph{simple splitting disk} for $\mu$ if $D$ satisfies the following conditions.
\begin{enumerate}
\item Each $I$--fiber of $N(B)$ intersects $D$ in at most one point.
\item $D\cap\mu=\emptyset$.
\item For any $I$--fiber $K$ that intersects $D$, both components of $K-D$ intersect $\mu$.
\end{enumerate}
Suppose $D$ is a simple splitting disk.  Let $N(B_\mu')$ and $N(B')$ be the manifold obtained by eliminating a small neighborhood of $D$ from $N(B_\mu)$ and $N(B)$ respectively.  So, we may consider $N(B_\mu')$ and $N(B')$ as fibered neighborhoods of branched surfaces $B_\mu'$ and $B'$ respectively. $B_\mu'$ and $B'$ are called the branched surfaces obtained by splitting along $D$.  By our assumptions on $D$, $B_\mu'$ is the sub-branched surface of $B'$ that fully carries $\mu$.  It is possible that some surfaces in $\{S_n\}$ are not carried by $B'$ anymore.  Nonetheless, we have the following Lemma.  

Recall that if $\mu\subset N(B)$ is a lamination carried by $B$, and $B'$ is obtained by splitting $B$, then we may assume $N(B')\subset N(B)$ and we say that $\mu$ is carried by $B'$ if $\mu\subset N(B')$ after some $B$--isotopy (see section~\ref{Spre} for the definition of $B$--isotopy).

\begin{lemma}\label{Lsdisk}
Let $B$, $\mu$, $B'$, $D$ and $\{S_n\}$ be as above.
There are a finite collection of branched surfaces, obtained by splitting $B$, such that
\begin{enumerate}
\item each $S_n$ is carried by a branched surface in this collection,
\item $B'$ is in the collection,
\item if another branched surface $B''$ in this collection carries $\mu$, then $B'$ is a sub-branched surface of $B''$.  In particular, $B'$ and $B''$ have the same sub-branched surface that fully carries $\mu$.
\end{enumerate}
\end{lemma}
\begin{proof}
Let $E$ be the union of $I$--fibers of $N(B)$ that intersect $D$. So, $E=\pi^{-1}(\pi(D))$, where $\pi:N(B)\to B$ is the collapsing map.  Since each $I$--fiber of $N(B)$ intersects $D$ in at most one point, $E$ is homeomorphic to a 3--ball $D^2\times I$.  After some small perturbation, we may simply identify $E$ to $D^2\times I$ with each $I$--fiber of $E$ coming from an $I$--fiber of $N(B)$.

If $B'$ carries every surface in $\{S_n\}$, then there is nothing to prove.  Suppose $S_n$ is not carried by $B'$.  Then $S_n\cap D\ne\emptyset$ under any $B$--isotopy.  If $S_n\cap\mu=\emptyset$ in $N(B)$, then by adding some branch sectors to $B'$, we can construct a branched surface $B''$ that carries $S_n$, and $B''$ satisfies part 3 of the lemma (this construction is similar to Figure~\ref{sp1}, where one can obtain $\tau_1$ by adding a branch sector to $\tau_2$).  Moreover, $B''$ can also be obtained by splitting $B$.  Since $D$ is compact, there are only finitely many ways to add such branch sectors.  Hence, there are only finitely many such $B''$.

Next, we will assume $S_n\cap\mu\ne\emptyset$ under any $B$--isotopy.   $S_n\cap E$ is a union of compact surfaces transverse to the $I$--fibers and each component of $S_n\cap E$ is $B$--isotopic to a sub-surface of $D$.  Let $P$ be a component of $S_n\cap E$ such that $P\cap\mu\ne\emptyset$ under any $B$--isotopy.  We may assume $P$ intersects both components of $D^2\times\partial I$, where $D^2\times I=E$ as above.  Since $S_n\cap\mu\ne\emptyset$ under any $B$--isotopy, there is a relatively short arc $\alpha\subset P$ with endpoints in different components of $D^2\times\partial I$, and after slightly extending $\alpha$ in $S_n$, we may assume $\alpha\cap\mu\ne\emptyset$ under any $B$--isotopy.  So, by deleting a small neighborhood of $\alpha$ from $N(B)$, we can split $B$ into a new branched surface $B_1$.  This splitting is similar to the splitting from $\tau$ to $\tau_1$ in Figure~\ref{sp1}.  By the construction, $B_1$ carries $S_n$, but since $\alpha\cap\mu\ne\emptyset$ under any $B$--isotopy, $B_1$ does not carry $\mu$.  Since $D$ is fixed, up to $B$--isotopy, there are only finitely many such compact surfaces $S_n\cap E$, and there are only finitely many different splittings like this.  Hence, we can perform such splittings on $B$ in a neighborhood of $E$ and obtain finitely many branched surfaces $B_1,\dots, B_k$, such that no $B_i$ carries $\mu$.  These $B_i$'s plus the branched surfaces $B''$ above are the collection of branched surfaces satisfying the conditions in the lemma.
\end{proof}

Note that any splitting along $\mu$ can be decomposed as a sequence of successive splittings along simple splitting disks.  Hence, we can apply Lemma~\ref{Lsdisk} at each step and obtain a collection of branched surfaces with similar properties.

\begin{lemma}\label{Ltorus0}
Let $B$ be a branched surface in $M$, and $T$ a compact orientable surface carried by $N(B)$.  Suppose $T$ is either a closed surface or a surface whose boundary lies in $\partial_vN(B)$.  Then, there is a finite collection of branched surfaces $B_1,\dots, B_k$ obtained by splitting $B$, such that 
\begin{enumerate}
\item if $B_i$ still carries $T$, then each $I$--fiber of $N(B_i)$ intersects $T$ in at most one point,
\item any closed surface carried by $B$ is carried by some $B_i$. 
\end{enumerate}
\end{lemma}
\begin{proof}
If every $I$--fiber of $N(B)$ intersects $T$ in at most one point, then there is nothing to prove.  Let $m$ ($m>1$) be the maximal number of points that an $I$--fiber of $N(B)$ intersects $T$, and let $I_m$ be the union of those $I$--fibers of $N(B)$ that intersect $T$ in $m$ points.  Since $m$ is maximal, $I_m$ is an $I$--bundle over a compact surface $F_m\subset N(B)$ and each $I$--fiber of $N(B)$ intersects $F_m$ in at most one point.  After ``blowing air" into $N(B)$ if necessary, we may assume $F_m$ is not a closed surface.  Moreover, since $m$ is maximal, $\partial F_m\subset \pi^{-1}(L)$, where $L$ is the branch locus of $B$ and $\pi:N(B)\to B$ is the collapsing map, and the induced branch direction at $\partial F_m$ points into $F_m$.

Note that if $F_m$ is non-orientable, $I_m$ is a twisted $I$--bundle over $F_m$.  Since both $T$ and $M$ are orientable, no matter whether $F_m$ is orientable or not, we may assume  that $F_m\cap T=\emptyset$ and for any $I$--fiber $K$ that intersects $F_m$, both components of $K-F_m$ intersect $T$.  Suppose a component of $F_m$ is not a disk, then let $\alpha$ be a properly embedded essential arc in $F_m$.  We can split $B$ in a small neighborhood of $\alpha$, as described in section~\ref{Spre} and shown in Figure~\ref{sp1}, and obtain a finite collection of branched surfaces with the following properties. 
\begin{enumerate}
\item Any closed surface carried by $B$ is still carried by a branched surface in this collection. 
\item Suppose $T$ is carried by a branched surface $B'$ in this collection. Let $I_m'$ be the union of $I$--fibers of $N(B')$ that intersect $T$ in $m$ points, hence $I_m'$ is an $I$--bundle over a compact surface $F_m'$.  Then, $F_m'$ is homeomorphic to the surface obtained by cutting $F_m$ open along $\alpha$. 
\end{enumerate} 

Thus, after a finite number of splittings, we may assume $F_m$ is a collection of disks.  Now, similar to the splittings in the proof of Lemma~\ref{Lsdisk}, we can split the branched  surface in a neighborhood of each disk component of $F_m$.  Since each $I$--fiber intersects $F_m$ in at most one point, such splittings take place in disjoint 3--balls.  So, after these splittings, we get a collection of branched surfaces that satisfy part 2 of this lemma, and if $T$ is still carried by a branched surface $B_i$ in this collection, then the maximal number of points that an $I$--fiber of $N(B_i)$ intersects $T$ is smaller than $m$.  Therefore, we can apply these splittings to each branched surface in this collection, and eventually get $m=1$ for each branched surface that carries $T$.   
\end{proof}

Now, we consider a torus $T$ carried by a branched surface $B$ in $M$. 

\begin{lemma}\label{Ltorus1}
Let $B$ be a branched surface in $M$ and $T\subset N(B)$ an embedded torus carried by $N(B)$.  Suppose each $I$--fiber intersects $T$ in at most one point and $T$ bounds  a solid torus in $M$.  Let $S\subset N(B)$ be a closed orientable surface fully carried by $B$ and $S\cap T\ne\emptyset$ under any $B$--isotopy.  Then, there are a surface $S'$, a number $\sigma$, and an arc $\alpha\subset S'$ such that
\begin{enumerate}
\item $S'$ is carried by $B$ and is isotopic to $S$ in $M$,
\item $length(\alpha)<\sigma$ and $\sigma$ depends only on $B$ and $T$, not on $S$.
\item $\alpha\cap T\ne\emptyset$ under any $B$--isotopy, 
\end{enumerate} 
\end{lemma}
\begin{proof}
Let $E$ be the union of the $I$--fibers of $N(B)$ that intersect $T$.  Since each $I$--fiber intersects $T$ in at most one point, $E$ is homeomorphic to an $I$--bundle $T^2\times I$.  After some perturbation at $T^2\times\partial I$, we may assume the $I$--fibers of $E=T^2\times I$ are from the $I$--fibers of $N(B)$ and $T=T^2\times\{1/2\}\subset E$.  $S\cap E$ is a union of compact orientable surfaces properly embedded in $E$ and transverse to the $I$--fibers.  

Let $T_0$ and $T_1$ be the two components of $T^2\times\partial I$.  If the boundary of every component of $S\cap E$ lies in the same component of $T^2\times\partial I$, then after some $B$--isotopy, $T$ is disjoint from $S$, which contradicts our hypothesis.  So, there must be a component of $S\cap E$, say $P$, intersecting both $T_0$ and $T_1$.  

Suppose a component $c$ of $\partial P$ is a trivial circle in $T_i$ and let $\Delta_c$ be the disk in $T_i$ bounded by $c$.  We say the circle $c$ is of type $I$, if (after smoothing out the corner) $P\cup\Delta_c$ is a surface transverse to the $I$--fibers of $E$, otherwise, $c$ is of type $II$.  For each innermost trivial circle $c$ in $\partial P\cap T_i$ of type $I$, we can glue the disk $\Delta_c$ to $P$ and then push (a neighborhood of) the disk into the interior of $T^2\times I$.  This operation yields a new surface transverse to the $I$--fibers of $T^2\times I$. We can keep performing such operations on the resulting surface and eventually get a surface $\hat{P}$ such that $\partial\hat{P}$ contains no trivial circle of type $I$.  $\hat{P}$ is a connected compact surface properly embedded in $T^2\times I$ and transverse to the $I$--fibers.  We have the following 4 cases to consider.

Case 1. $\partial\hat{P}$ contains a trivial circle in $T_i$.

Let $c$ be an innermost trivial circle of $\partial\hat{P}\cap T_i$, and $c$ bounds a disk $\Delta_c$ in $T_i$. By the assumptions on $\hat{P}$, $c$ is of type $II$. Now, we cut $E=T^2\times I$ open along $\hat{P}$ and obtain a manifold $N$ which is the closure (under path metric) of $E-\hat{P}$.  Since $\hat{P}$ is transverse to the $I$--fibers, we may consider $N$ as an induced $I$--bundle with its vertical boundary pinched into circles/cusps.  Let $N_1$ be the component of $N$ containing $\Delta_c$.  Since $c$ is of type $II$, $\Delta_c$ must be a component of the horizontal boundary of the pinched $I$--bundle $N_1$.  Thus, $N_1$ is a product $D^2\times I$ with vertical boundary $\partial D^2\times I$ pinched to a circle.  As $\hat{P}$ is connected, $\hat{P}$ must be a disk $B$--isotopic to $\Delta_c\subset T_i$.  

Since $\hat{P}$ is obtained by gluing disks to $P$, $P$ must be a planar surface $B$--isotopic to a sub-surface of $T$.  Moreover, by our assumptions on $P$, $\partial P$ has components in both $T_0$ and $T_1$.  Thus, there is an arc $\alpha$ properly embedded in $P$ connecting a component of $\partial P\cap T_0$ to a component of $\partial P\cap T_1$.  Since $P$ is $B$--isotopic to a sub-surface of $T$, we can choose $\alpha$ so that $length(\alpha)$ is bounded from above by a number $\sigma$ that depends only on $T$ and $B$, not on $S$.

Case 2. $\partial\hat{P}=\emptyset$.

Since $\hat{P}$ is transverse to $I$--fibers, this implies that $\hat{P}$ is a torus $B$--isotopic to $T$.  Since $\hat{P}$ is obtained by gluing disks to $P$, $P$ is $B$--isotopic to a sub-surface of $T$.  Since $\partial P$ has components in both $T_0$ and $T_1$, as in case 1, we can find an arc $\alpha\subset P$ connecting a component of $\partial P\cap T_0$ to a component of $\partial P\cap T_1$, and the length of $\alpha$ is bounded by a number $\sigma$ that does not depend on $S$.

Case 3. $\partial\hat{P}$ contains no trivial circle and $\partial\hat{P}\subset T_0$.

In this case, $\partial\hat{P}$ consists of parallel essential simple closed curves in the torus $T_0$.  As in case 1, we cut $E=T^2\times I$ open along $\hat{P}$ and obtained a pinched $I$--bundle $N$. Since $\partial\hat{P}\subset T_0$, $N$ has a component $N_1$ containing $T_1$.  As $\partial\hat{P}\subset T_0$, $T_1$ is a component of the horizontal boundary of $N_1$.  Hence, $\hat{P}$ is $B$--isotopic to a sub-surface of $T_1$.  Moreover, since $\hat{P}$ is connected and $\partial\hat{P}$ contains no trivial circle, $\hat{P}$ is an annulus with $\partial\hat{P}\subset T_0$.  Since $\hat{P}$ is obtained by gluing disks to $P$, $P$ is a planar surface $B$--isotopic to a sub-surface of $T$.  As in case 1, we can find an arc $\alpha\subset P$ connecting a component of $\partial P\cap T_0$ to a component of $\partial P\cap T_1$, and the length of $\alpha$ is bounded by a number $\sigma$ that does not depend on $S$.

Case 4. $\partial\hat{P}$ contains no trivial circle and $\partial\hat{P}$ has components in both $T_0$ and $T_1$.

As before, let $N$ be the manifold obtained by cutting $E=T^2\times I$ open along $\hat{P}$, and $N$ is a pinched $I$--bundle with the bundle structure induced from that of $T^2\times I$.  The two sides of $\hat{P}$ correspond to two sub-surfaces $\hat{P}^+$ and $\hat{P}^-$ in the horizontal boundary of $N$.  As $\partial\hat{P}$ contains no trivial circle, $\partial\hat{P}^+$ and $\partial\hat{P}^-$ does not bound disks in the horizontal boundary of $N$.  So, $\hat{P}^+$ and $\hat{P}^-$ are $\pi_1$--injective in $N$, which implies that $\hat{P}$ is incompressible and $\pi_1$--injective in $E=T^2\times I$.  So, $\pi_1(\hat{P})$ is a subgroup of $\mathbb{Z}\oplus\mathbb{Z}$.  By the assumption on $\hat{P}$ in this case, $\hat{P}$ must be an annulus with one boundary circle in $T_0$ and the other boundary circle in $T_1$.

If the distance (in $\hat{P}$) between the two components of $\partial\hat{P}$ is large, then the annulus $\hat{P}$ wraps around $T$ many times.  Since $\hat{P}$ is obtained by gluing disks to $P$, either there is a relatively short arc $\alpha$ properly embedded in $P$ connecting $\partial P\cap T_0$ to $\partial P\cap T_1$, or $P$ contains a sub-surface which is a long annulus wrapping around $T$ many times.  In the latter case, we can perform a Dehn twist in $T^2\times I$ to unwrap $\hat{P}$ and $P$.  Since $T$ bounds a solid torus in $M$, a Dehn twist around $T$ is an isotopy in $M$.  Therefore, after a Dehn twist in $T^2\times I$, we get a surface $S'$, which is isotopic to $S$ in $M$ and also fully carried by $B$, such that there is an arc $\alpha$ connecting $S'\cap T_0$ to $S'\cap T_1$ and $length(\alpha)$ is less than a fixed number $\sigma$ that does not depend on $S$ or $S'$.  

After slightly extending such arcs $\alpha$ in $S$ or $S'$, we have $\alpha\cap T\ne\emptyset$ under any $B$--isotopy.
\end{proof}

\begin{lemma}\label{Ltorus2} 
Let $B$ be a branched surface in $M$, $T\subset N(B)$ an embedded torus carried by $N(B)$, and suppose $T$ bounds  a solid torus in $M$.  Let $\{S_n\}$ be a sequence of closed orientable surfaces carried by $B$ and with genus at least 2.  Then, there is a finite collection of branched surfaces, obtained by splitting $B$ and then taking sub-branched surfaces, with the following properties. 
\begin{enumerate}
\item No branched surface in this collection carries $T$. 
\item For each $S_n$, there is a surface $S_n'$ isotopic to $S_n$ in $M$ and fully carried by a branched surface in this collection.
\end{enumerate}
\end{lemma}
\begin{proof}
This lemma is an easy corollary of Lemmas \ref{Ltorus0} and \ref{Ltorus1}. 
If there is an $I$--fiber of $N(B)$ that intersects $T$ in more than one point, by Lemma~\ref{Ltorus0}, we can split $B$ into a finite collection of branched surfaces $B_1,\dots, B_m$, such that any surface carried by $B$ is carried by some $B_i$, and if $B_i$ carries $T$, each $I$--fiber of $N(B_i)$ intersects $T$ in at most one point.  Moreover, after taking sub-branched surfaces of each $B_i$, we may also assume that each $S_n$ is fully carried by some $B_i$.

First note that if a branched surface $B$ fully carries $S_n$ then no component of $\partial_hN(B)$ can be a torus.  This is because if $B$ fully carries $S_n$, $S_n$ intersects every $I$--fiber of $N(B)$ and hence any component of $\partial_hN(B)$ is isotopic to a subsurface of $S_n$.  Thus if $\partial_hN(B)$ has a torus component then $S_n$ must be a torus, contradicting that $S_n$ has genus at least two.

Let $B_i$ be a branched surface in this collection that carries $T$ and fully carries $S_n$.  If $S_n\cap T=\emptyset$ in $N(B_i)$, then we cut $N(B_i)$ open along $T$ and obtain $\overline{N(B_i)-T}$ which carries both $T$ and $S_n$.  However, a horizontal boundary component of $\overline{N(B_i)-T}$ is a torus parallel to $T$ and the argument above implies that $\overline{N(B_i)-T}$ dose not fully carry $S_n$.  So, after taking a sub-branched surface of $\pi(\overline{N(B_i)-T})$ (where $\pi$ is the map collapsing every $I$--fiber of $\overline{N(B_i)-T}$ to a point), we get a branched surface $B_i'$ that fully carries $S_n$.  Note that the operation of taking a sub-branched surface destroys the torus components of the horizontal boundary that come from cutting $N(B_i)$ along $T$. So the new branched surface $B_i'$ does no carry $T$.  A branched surface has only finitely many sub-branched surfaces.  Thus, after these operations, we may assume each $B_i$ has the property that $S_n\cap T\ne\emptyset$ under any $B_i$--isotopy, if $B_i$ carries $T$ and fully carries $S_n$.  Now, by Lemma~\ref{Ltorus1}, for each surface $S_n$ fully carried by $B_i$, we can find a surface $S_n'$ and an arc $\alpha\subset S_n'$, such that $S_n'$ is isotopic to $S_n$ in $M$, $S_n'$ is also fully carried by $B$, $\alpha\cap T\ne\emptyset$ under any $B_i$--isotopy, and $length(\alpha)$ is bounded from above by a fixed number $\sigma$ depending only on $T$ and $B_i$.    

Then, similar to the proof of Lemma~\ref{Lsdisk}, we split $N(B_i)$ in a small neighborhood of $\alpha$, as the splitting from $\tau$ to $\tau_1$ in Figure~\ref{sp1}.  Since $\alpha\subset S_n'$ and $\alpha\cap T\ne\emptyset$ under any $B_i$--isotopy, we may perform the splitting so that the branched surface after this splitting still carries $S_n'$ but does not carry $T$.

We may assume $\pi(\alpha)$ is transverse to the branch locus.  Since $T$ is fixed and the length of $\alpha$ is bounded by a number $\sigma$ which depends only on $B_i$ and $T$, there are only a finite number of different such splittings along arcs like $\alpha$.  Thus, after performing a finite number of splittings on $B_i$, we get a finite collection of branched surfaces with the following properties.
\begin{enumerate}
\item No branched surface in this collection carries $T$. 
\item For any surface $S_n$ carried by $B_i$, there is a surface $S_n'$ that is isotopic to $S_n$ in $M$ and carried by $B_i$. 
\end{enumerate}

After performing these splittings on each $B_i$ and taking sub-branched surfaces if necessary, we get a collection of branched surfaces satisfying the properties in the lemma.
\end{proof}

\section{Splitting branched surfaces, the lamination case}\label{Slam}

Suppose $M$ is a closed, orientable, irreducible and atoroidal 3--manifold, and $M$ is not a Seifert fiber space.  By \cite{JR}, we may assume $M$ has a 0--efficient triangulation.  By \cite{R,St}, every strongly irreducible Heegaard surface is isotopic to an almost normal surface with respect to the 0--efficient triangulation.   As in section~\ref{Spre} and Proposition~\ref{Pfinite}, we can construct a finite collection of branched surfaces by gluing together normal disks and almost normal pieces, and each strongly irreducible Heegaard surface is fully carried by a branched surface in this collection.  Since the only normal 2--sphere is the vertex-linking one, after taking sub-branched surfaces, we may assume no branched surface in this collection carries any normal 2--sphere.  

Let $B$ be a branched surface in this collection.  Since $B$ fully carries an almost normal surface, at most one branch sector of $B$ contains an almost normal piece.  Let $b$ be the branched sector of $B$ that contains the almost normal piece.  As in section~\ref{Spre}, $B_N=B-int(b)$ is a sub-branched surface of $B$, which is called the normal part of $B$.  Clearly, every surface carried by $B_N$ is a normal surface, and $B_N$ does not carry any 2--sphere.   

\begin{lemma}\label{Lmu}
Let $B$ and $B_N$ be the branched surface constructed above, and let $\{S_n\}$ be a sequence of strongly irreducible Heegaard surfaces fully carried by $B$.  Let $\mu$ be an exceptional minimal measured lamination carried by $B_N$ with $\chi(\mu)=0$.  Then, $B$ can be split into a finite collection of branched surfaces with the following properties.
\begin{enumerate}
\item Up to isotopy, each $S_n$ is carried by a branched surface in this collection.
\item No branched surface in this collection carries $\mu$.
\end{enumerate}
\end{lemma}
\begin{proof}
Let $B_\mu$ be the sub-branched surface of $B$ that fully carries $\mu$.  So, $B_\mu$ is also a sub-branched surface of $B_N$.  Since $B_N$ does not carry any 2--sphere, $B_\mu$ does not carry any 2--sphere either.  Moreover, every torus $T$ carried by $B_\mu$ is a normal torus, and by Lemma~\ref{Ltorus}, $T$ bounds a solid torus in $M$.

Next, we perform some splittings on $B_\mu$.  By Lemma~\ref{Llocus}, after some splittings, we have the following.
\begin{description}
\item[Property A] the branch locus $L_\mu$ of $B_\mu$ has no double point, 
\item[Property B] $B_\mu-L_\mu$ consists of annuli and M\"{o}bius bands.
\end{description}  
Note that any splitting on a branched surface can be divided into a sequence of successive splittings along simple splitting disks (see section~\ref{Storus} for  definition).   By Lemma~\ref{Lsdisk}, we can perform splittings on $B$ and $B_\mu$ and obtain a finite collection of branched surfaces, such that 
\begin{enumerate}
\item each surface in $\{S_n\}$ is carried by a branched surface in this collection,
\item if a branched surface $B'$ in this collection carries $\mu$, then $B_\mu$, the sub-branched surface of $B'$ fully carrying $\mu$, satisfies properties A and B above.
\end{enumerate}

To simplify notation, we still use $B$ to denote a branched surface in this collection that carries $\mu$, use $B_\mu$ to denote the sub-branched surface of $B$ fully carrying $\mu$, and assume $B_\mu$ satisfies properties A and B above.  Moreover, as in the proof of Lemma~\ref{Llocus}, we may assume the number of components of $M-N(B_\mu)$ is minimal among branched surfaces fully carrying $\mu$.  After some isotopy, we can also assume $\partial_hN(B_\mu)\subset\mu$.

Since $B_\mu$ is a sub-branched surface of $B$, we may also consider $B$ as a branched surface obtained by adding some branch sectors to $B_\mu$.  Next, we will fix $B_\mu$ and split $B$ near $B_\mu$.   

We first analyze how the branch sectors of $\overline{B-B_\mu}$ are added to $B_\mu$ at the cusps of $B_\mu$.  Let $c_x$ be a circle in $\partial_hN(B_\mu)$ parallel and close to a boundary circle of $\partial_hN(B_\mu)$.   Since the branch locus of $B_\mu$ contains no double point, $\pi(c_x)$ is a circle parallel and close to a component of the branch locus of $B_\mu$.  To simplify notation, we use $l_x$ to denote both the component of the branch locus and the corresponding cusp.  The union of the $I$--fibers of $N(B_\mu)$ that intersect $c_x$ is a vertical annulus $A_x$, and $A_x\cap\mu$ is a union of parallel circles.  By assuming $N(B_\mu)\subset N(B)$ as before, we may consider $A_x$ as a vertical annulus in $N(B)$.  Let $\hat{A}_x$ be the union of $I$--fibers of $N(B)$ that intersect $A_x$.  After enlarging $\hat{A_x}$ a little, we may consider $\hat{A}_x$ as a fibered neighborhood of a train track $\tau_x$ which consists of the circle $\pi(c_x)$ and some ``tails" along the circle, as the top train track in Figure~\ref{sp2}.  Note that $\tau_x$ can be regard as the ``spine" of a small neighborhood of $\pi(c_x)$ in $B$.  If $S_n$ is fully carried by $B$, then $S_n\cap\hat{A}_x$ is a union of arcs and/or circles transverse to the $I$--fibers.  We have 3 cases.   

Case 1 is that $S_n\cap\hat{A}_x$ contains a circle.  Case 2 is that $S_n\cap\hat{A}_x$ contains a spiral wrapping around $\hat{A}_x$ more than twice.  Case 3 is that $S_n\cap\hat{A}_x$ contains no circle and the length of every arc in $S_n\cap\hat{A}_x$ is relatively short (compared with the length of the circle $\pi(c_x)$ in the train track).  

Now, we split $\hat{A}_x$ along $S_n\cap\hat{A}_x$.  In the first case, we can split $\hat{A}_x$ along some relatively short arcs, as shown in splitting 1 of Figure~\ref{sp2}, and get a vertical annulus whose intersection with $S_n$ consists of circles.  In the second case, we can split $\hat{A}_x$ along relatively short arcs, as shown in splitting 2 of Figure~\ref{sp2}, and get a fibered neighborhood of a train track whose intersection with $S_n$ consists of only spirals.  The train track in the second case consists of a circle and some ``tails", and on each side of the circle, the cusps of the ``tails" have the same direction.  In the third case, as shown in splitting 3 of Figure~\ref{sp2}, the splitting along a short arc will destroy the annulus $A_x$ and the circle $\pi(c_x)=\pi(A_x)$.

\begin{figure}
\begin{center}
\psfrag{1}{1}
\psfrag{2}{2}
\psfrag{3}{3}
\includegraphics[width=4in]{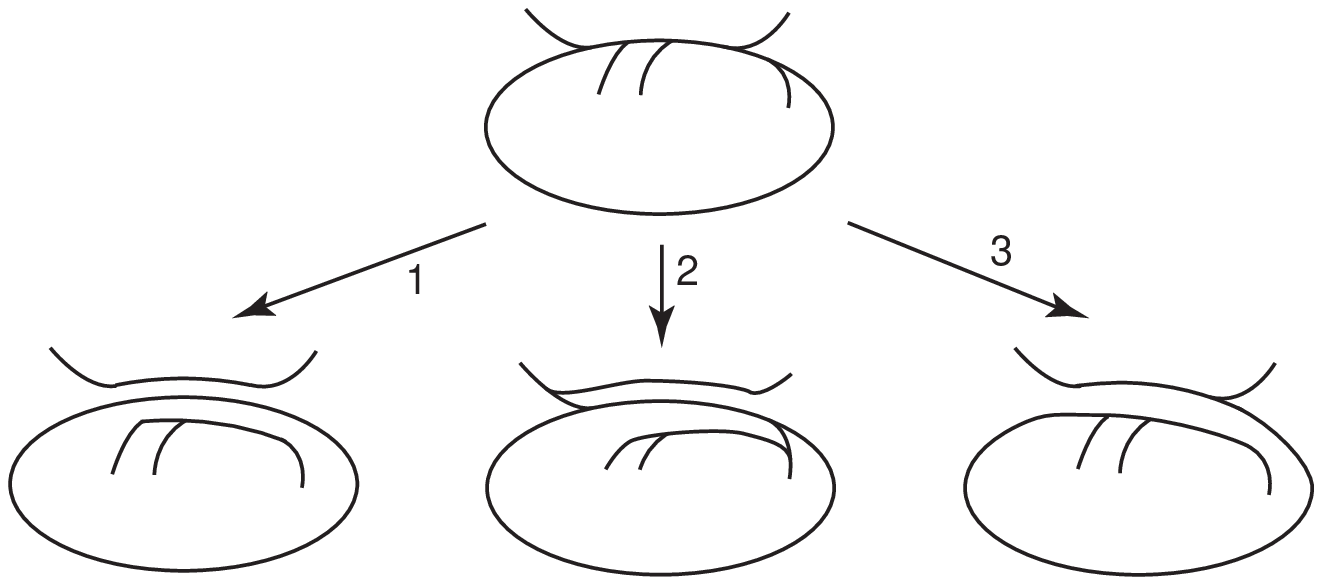}
\caption{} \label{sp2}
\end{center}
\end{figure}

Now, we consider the bigger 3-dimensional pictures of the splittings above. The third case is simple. In the third case, similar to the splittings in the poof of Lemma~\ref{Lsdisk}, the branched surface after the splitting does not carry $\mu$ anymore. Next, we will focus on the first two cases.

For the first two cases, since $c_x$ lies in a small neighborhood of a boundary circle of $\partial_hN(B_\mu)$, we may assume the (2-dimensional) splittings occur in a small neighborhood of the cusp, and we perform some splittings and pinchings on $B$ accordingly, as shown in Figure~\ref{cusp} (a).  Although both local splittings in Figure~\ref{cusp} (a) may happen in the two cases, the basic picture for the splittings in case 1 is the splitting 1 in Figure~\ref{cusp} (a), and the basic picture for the splittings in case 2 is the splitting 2 in Figure~\ref{cusp} (a).  

To simplify notation, we still use $B$ to denote the branched surface after the splittings above.  In the first two cases, $B_\mu$ is still a sub-branched surface of $B$ after the splittings.  In case 1, no branch sector of $\overline{B-B_\mu}$ intersects the cusp $l_x$ after the splittings, in other words, the cusp $l_x$ for $B_\mu$ is a cusp for $B$ after the splitting.  In case 2, as shown in Figure~\ref{cusp} (b), the branch sectors of $\overline{B-B_\mu}$ have the coherent direction along the cusp $l_x$ after the splittings.  

\begin{figure}
\begin{center}
\psfrag{1}{1}
\psfrag{2}{2}
\psfrag{g}{$\gamma_c$}
\psfrag{(a)}{(a)}
\psfrag{(b)}{(b)}
\includegraphics[width=4in]{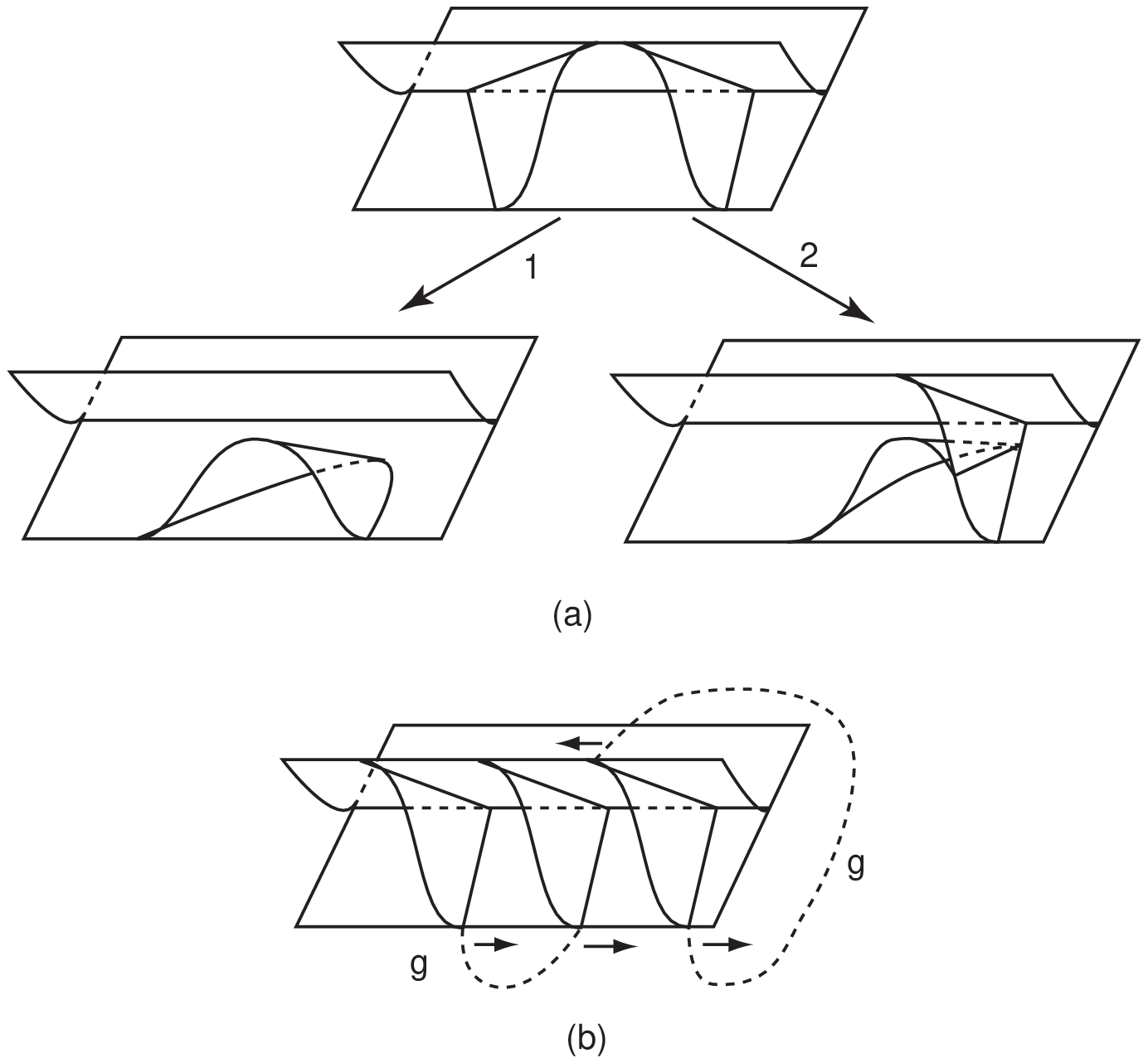}
\caption{} \label{cusp}
\end{center}
\end{figure}

Next, we show that the second case cannot happen.  To prove this, we first show that the second case cannot happen at the cusp of a $monogon\times S^1$ region of $M-int(N(B_\mu))$.  By Lemma~\ref{L0eff}, there is always a $monogon\times S^1$ region in $M-int(N(B_\mu))$.  Let $D$ be a monogon disk properly embedded in $M-int(N(B_\mu))$ with $\partial D=\alpha\cup\beta$, where  $\beta\subset\partial_vN(B_\mu)$ is a vertical arc in $\partial_vN(B_\mu)$ and $\alpha\subset\partial_hN(B_\mu)$.  We can use $D\times S^1$ to denote the $monogon\times S^1$ region of $M-int(N(B_\mu))$.  So, $\beta\times S^1$ is a component of $\partial_vN(B_\mu)$.  As $\partial_hN(B_\mu)\subset\mu$, $\alpha\times S^1$ lies in a leaf $l$ of $\mu$.  As in the proof of Lemma~\ref{Llocus}, $l$ is an infinite annulus.  Since the number of components of $M-N(B_\mu)$ is minimal, $l$ is the boundary (under path metric) of a component of $M-\mu$ which is the product of $S^1$ and an end-compressing disk (i.e. a monogon with an infinitely long tail, see page 45 of \cite{GO}).  

Suppose we are in case 2 at the cusp of the $monogon\times S^1$ region $D\times S^1$, and suppose the branch sectors of $\overline{B-B_\mu}$ are coherent along this cusp $\beta\times S^1$, as shown in Figure~\ref{cusp} (b).  For any surface $S_n$ fully carried by $B$, we can regard $S_n\cap (D\times S^1)$ as a compact surface carried by those branch sectors of $\overline{B-B_\mu}$ in this $monogon\times S^1$ region.  Let $C_n$ be a component of $S_n\cap (D\times S^1)$ whose boundary intersects the cusp $\beta\times S^1$.  The union of $C_n$ and $B_\mu$ naturally form a sub-branched surface of $B$.  Let $c$ be a boundary circle of $C_n$ that intersects the cusp of this $monogon\times S^1$ region.  So, $c$ is a circle lying in the branch locus of $B$ and has an induced cusp/branch direction.  Let $\gamma_c$ be an arc in $c$ with both endpoints in the cusp $l_x=\pi(\beta\times S^1)$, see the two dashed arcs in Figure~\ref{cusp}~(b) for pictures of $\gamma_c$.  We may regard $\gamma_c$ as an arc properly embedded in the annulus $\alpha\times S^1\subset\partial_hN(B_\mu)$ and $\gamma_c$ has a normal direction induced from the cusp direction at $\gamma_c$.  Since we are in case 2 and the branch sectors of $\overline{B-B_\mu}$ are coherent along this cusp $l_x$, as shown in Figure~\ref{cusp}~(b), the cusp directions at $\partial\gamma_c$ cannot be extended to a compatible cusp direction along $\gamma_c$.  Hence, the second case can never happen near the cusp of a $monogon\times S^1$ region.  In other words, after some splittings as in Figure~\ref{cusp}~(a), either $\overline{B-B_\mu}$ has no branch sector intersecting the cusp of any $monogon\times S^1$ region of $M-B_\mu$, or the branched surface after the splitting does not carry $\mu$ anymore.

Now, we consider the cusp of any component $l_x$ of the branch locus of $B_\mu$, and suppose we are in case 2 at this cusp. So, we may assume the branch sectors of $\overline{B-B_\mu}$ at the cusp have coherent directions as in Figure~\ref{cusp} (b). Let $A_x'=\pi^{-1}(l_x)$, where $\pi: N(B_\mu)\to B_\mu$ is the collapsing map.  Since the branch locus of $B_\mu$ has no double points and $\mu$ is a measured lamination, $A_x'\cap\mu$ is a union of circles.  Since every leaf of $\mu$ is dense, as in the proof of Lemma~\ref{Llocus}, there is an annulus in $N(B_\mu)-\mu$, transverse to the $I$--fibers and connecting $\beta\times S^1$ to $A_x'$, where $\beta\times S^1$ is the cusp of a $monogon\times S^1$ region above.  By deleting a small neighborhood of this annulus, we can split $B_\mu$ so that the cusp of this $monogon\times S^1$ region passes $l_x$ and lies in a small neighborhood of $l_x$, as shown in Figure~\ref{sp3} (a).  Since no branch sector of $\overline{B-B_\mu}$ intersects the cusp of a $monogon\times S^1$ region, this splitting does not really affect the branch sectors of $\overline{B-B_\mu}$.

\begin{figure}
\begin{center}
\psfrag{m}{$monogon\times S^1$}
\psfrag{s}{splitting}
\psfrag{l}{$l_x$}
\psfrag{c}{$c_x$}
\psfrag{(a)}{(a)}
\psfrag{(b)}{(b)}
\includegraphics[width=4in]{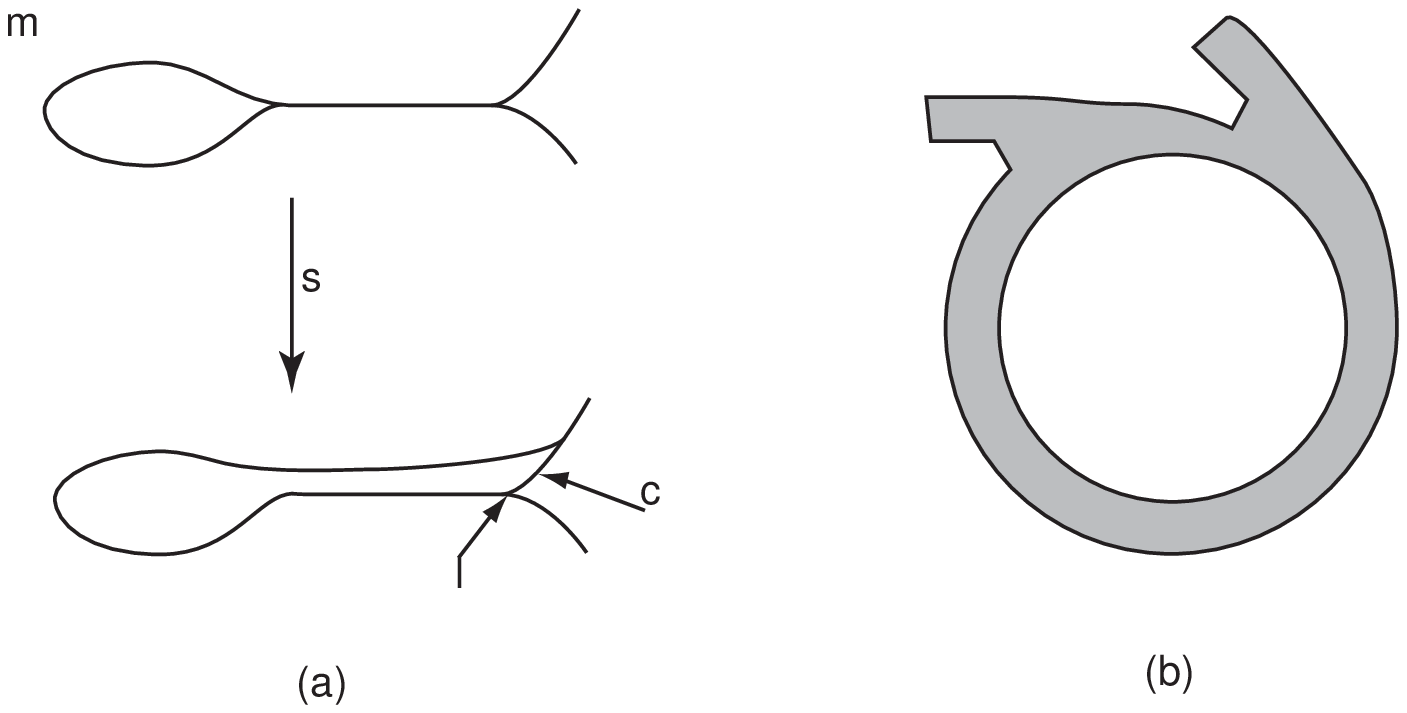}
\caption{} \label{sp3}
\end{center}
\end{figure}

As before, let $c_x$ be a circle in $\partial_hN(B_\mu)$ parallel and close to the cusp, and let $A_x=\pi^{-1}(c_x)$ be the vertical annulus, where $\pi:N(B_\mu)\to B_\mu$ is the collapsing map.  By our assumptions on the branch locus of $B_\mu$, $\partial A_x$ lies in $\partial_hN(B_\mu)$. Since we have split $B_\mu$ along the annulus above, as shown in Figure~\ref{sp3} (a), we may assume a component of $\partial A_x$ lies in the horizontal boundary of the $monogon\times S^1$ region in $M-int(N(B_\mu))$, close to the cusp $\beta\times S^1$.    

We may consider $N(B_\mu)\subset N(B)$ and consider $A_x$ as a vertical annulus in $N(B)$.  Let $\hat{A}_x$ be the union of $I$--fibers of $N(B)$ that intersect $A_x$.  As before, after enlarging $\hat{A}_x$ a little, we may consider $\hat{A}_x$ as a fibered neighborhood of a train track which consists of a circle and some ``tails" along the circle.  

Since we are in case 2, by our assumptions above, the branch sectors of $\overline{B-B_\mu}$ that intersect the cusp of $l_x$ have coherent direction, and there is no branch sector of $\overline{B-B_\mu}$ intersecting the cusp of the $monogon\times S^1$ region.  So, a neighborhood of $\hat{A}_x$ in $N(B)$ must be like Figure~\ref{sp3} (b), where the smooth boundary circle is a circle in $\partial_hN(B_\mu)$ parallel and close to the cusp of the $monogon\times S^1$ region and the ``tails" correspond to the branch sectors of $\overline{B-B_\mu}$ that intersect the cusp of $l_x$.  Since the branch sectors of $\overline{B-B_\mu}$ that intersect the cusp of $l_x$ have coherent direction, these ``tails" in $\hat{A}_x$ have the same cusp/branch direction along the annulus $A_x$, as shown in Figure~\ref{sp3}~(b).  A standard Poincar\'{e}-Bendixson type argument implies that any curves fully carried by $\hat{A}_x$ must contain an infinite spiral and a limit cycle.  So, $B$ cannot fully carry any compact surface with this configuration.  Therefore, the second case above cannot happen at any cusp circle of $B_\mu$.  

Recall that in case 3, after the splittings above, the branched surface does not carry $\mu$ anymore.  In case 1, after the splittings, $\overline{B-B_\mu}$ does not contain any branch sector that intersects the branch locus $L_\mu$ of $B_\mu$, which means each circle in $L_\mu$ is a component of the branch locus of $B$.

Since the splittings performed above are along relatively short arcs, similar to the proof of Lemmas~\ref{Lsdisk} and \ref{Ltorus2}, by performing a finite number of splittings on $B$ (and taking sub-branched surfaces if necessary), we can obtain a finite collection of branched surfaces with the following properties.   
\begin{enumerate}
\item Each $S_n$ is fully carried by a branched surface in this collection,
\item If a branched surface $B'$ in this collection carries $\mu$, then $B_\mu$, the sub-branched surface of $B'$ that fully carries $\mu$, satisfies properties A and B before, and no branch sector of $\overline{B'-B_\mu}$ intersects the branch locus of $B_\mu$.
\end{enumerate}

To simplify notation, we still use $B$ to denote a branched surface in this collection that carries $\mu$, and use $B_\mu$ to denote the sub-branched surface of $B$ that fully carries $\mu$.  Let $L_\mu$ be the branch locus of $B_\mu$.  By our assumptions before, $B_\mu-L_\mu$ consists of annuli and M\"{o}bius bands.  We only need to consider the case that $B_\mu-L_\mu$ is a union of annuli, and the proof for the case that $B_\mu-L_\mu$ contains a M\"{o}bius band is the same after ``blowing air" into $N(B_\mu)$.  Let $l$ be an essential simple closed curve in a component of $B_\mu-L_\mu$. Let $A=\pi^{-1}(l)$, where $\pi:N(B_\mu)\to B$ is the collapsing map.  So, $A$ is a vertical annulus with $\partial A\subset\partial_hN(B_\mu)$, and $A\cap\mu$ is a union of parallel circles.   

By assuming $N(B_\mu)\subset N(B)$, we may consider the annulus $A$ as a vertical annulus in $N(B)$, and we denote the union of the $I$--fibers of $N(B)$ that intersect $A$ by $\hat{A}$.  As before, after enlarging $\hat{A}$ a little, we may consider $\hat{A}$ as a fibered neighborhood of a train track which consists of the circle $l$ and some ``tails" along $l$.  Suppose $S_n$ is fully carried by $B$.  We now split this train track along $S_n$.  As before, we have 3 cases as shown in Figure~\ref{sp2}.  Similarly, in the third case, the branched surface after splitting along a short arc does not carry $\mu$ anymore.  Since every leaf is dense in $\mu$, we can find an annulus in $N(B_\mu)-\mu$ connecting any component of $\partial_vN(B_\mu)$ to $A$.  Thus, by the argument on $L_\mu$ above, case 2 cannot happen at the circle $l$.  Therefore, we can split the branched surface into a finite collection of branched surfaces, such that each $S_n$ is carried by a branched surface in this collection, and if a branched surface $B$ in this collection carries $\mu$, then no branch sector of $\overline{B-B_\mu}$ intersects $l$.

We can apply this argument to any set of essential simple closed curves $l_1,\dots,l_m$ in $B_\mu-L_\mu$. So after a finite number of splittings, we may assume that no branch sector of $\overline{B-B_\mu}$ intersects any $l_i$. Now, $\overline{B-B_\mu}$ is a branched surface with boundary and the boundary of $\overline{B-B_\mu}$ is a train track in $B_\mu-L_\mu-\cup_{i=1}^ml_i$.  Since every surface carried by $B_\mu$ is normal, we can find enough such circles $l_i$ so that, after the splittings above and eliminating disks of contact, the train track $\overline{B-B_\mu}\cap(B_\mu-L_\mu)$ does not carry any trivial circle in $B_\mu-L_\mu$.  Since we can assume case 2 never happens along any such circles, by choosing enough such circles $l_i$ and after some more splitting and pinching, the boundary of $\overline{B-B_\mu}$ becomes a union of disjoint essential simple closed curves in the annuli $B_\mu-L_\mu$.   Note that all the splittings above are along relatively short arcs and small disks.  Similar to Lemma~\ref{Lsdisk}, we can perform a finite number of different splittings on $B$ and obtain a finite collection of branched surfaces such that 
\begin{enumerate}
\item each $S_n$ is fully carried by a branched surface in this collection,
\item if a branched surface $B'$ in this collection carries $\mu$, then $B_\mu$, the sub-branched surface of $B'$ fully carrying $\mu$, satisfies properties A and B before, and the boundary train track of $\overline{B'-B_\mu}$ consists of essential simple closed curves in $B_\mu-L_\mu$.
\end{enumerate}

By Proposition~\ref{PMS}, we may assume $B_\mu$ and $N(B_\mu)$ satisfy the hypotheses of Lemma~\ref{L0eff}.  So, by Lemma~\ref{L0eff}, there is a torus component $\Gamma$ of $\partial N(B_\mu)$ bounding a solid torus in $M$ and the solid torus contains $N(B_\mu)$.  $\Gamma$ is a union of annulus components of $\partial_vN(B_\mu)$ and $\partial_hN(B_\mu)$.  Let $N$ be the solid torus bounded by $\Gamma$ ($N(B_\mu)\subset N$). Let $l$ be an essential simple closed curve in an annulus component of  $\Gamma\cap\partial_hN(B_\mu)$. Since $B_\mu$ satisfies properties A and B before, by applying part 1 of Lemma~\ref {L0eff} to the components of $N-int(N(B_\mu))$, it is easy to see that $l$ does not bound a meridian disk of $N$. 

Suppose $B$ carries $\mu$ and fully carries infinitely many $S_n$'s.  By our assumptions above, the boundary of $\overline{B-B_\mu}$ consists of essential circles in  $B_\mu-L_\mu$.  So, for each $S_n$ fully carried by $B$, we may assume that $S_n\cap\Gamma$ consists of parallel essential non-meridian curves, and for any such $S_n$, the slope of  $S_n\cap\Gamma$ is the same as the slope of $l$ above.  By Theorem~\ref{Tsch} (a theorem of Scharlemann \cite{S}), $S_n\cap N$ consists of $\partial$--parallel annuli and possibly one other component, obtained from one or two $\partial$--parallel annuli by attaching an unknotted tube along an arc parallel to an arc in $\Gamma-S_n$.  We call the latter kind of component in Scharlemann's theorem an \emph{exceptional component}. An exceptional component is either a twice punctured torus or a planar surface with 4 boundary circles.  Note that, for an exceptional component, if one fixes the annuli, then there is only one way to attach the tube, up to isotopy.

Each component of $S_n\cap N$ is carried by $B\cap N$.  Since the boundary of $\overline{B-B_\mu}$ consists of simple closed curves, after some small perturbation, we may assume $B$ is transverse to the torus $\Gamma$, $B\cap\Gamma$ consists of parallel essential non-meridian simple closed curves, $N(B)\cap\Gamma$ consists of vertical annuli, and $\mu\subset N(B_\mu)\subset N$.  

For each $\partial$--parallel annulus $A$ in the solid torus $N$, $A$ is isotopic (fixing $\partial A$) to an annulus $A'$ in $\Gamma$ and we call $A'$ the \emph{image} of $A$ in $\Gamma$.  Let $A_1$ and $A_2$ be two $\partial$--parallel annuli in $N$ with $\partial A_i\subset N(B)\cap\Gamma$ ($i=1,2$).  We say that $A_1$ and $A_2$ are equivalent if $A_1$ is isotopic to $A_2$ via an isotopy of $N$ fixing $\Gamma-N(B)$.  Thus, there are only finitely many equivalence classes for $\partial$--parallel annuli with boundary in $N(B)\cap\Gamma$.  Now, we consider the exceptional components as in Scharlemann's theorem above, and we say that two exceptional components (from two Heegaard surfaces) are equivalent if they are isotopic via an isotopy of $N$ fixing $\Gamma-N(B)$.  Since the isotopy class of an exceptional component only depends on the annuli where the tube is attached, there are only finitely many equivalence classes for the exceptional components.

Let $A_1$ and $A_2$ be properly embedded and disjoint annuli in $N$ carried by $N(B)\cap N$, and $A_i'$ be the image of $A_i$ in $\Gamma$.  The solid torus bounded by $A_i\cup A_i'$ must contain at least one component of $\partial_hN(B)\cap N$.  Moreover, if $A_1'$ and $A_2'$ are nested, say $A_1'\subset A_2'$, and if $A_1$ and $A_2$ are not $B$--isotopic, then the solid torus between $A_1$ and $A_2$, i.e. the solid torus bounded by $A_1\cup A_2\cup(A_2'-A_1')$, must contain at least one component of $\partial_hN(B)\cap N$.  Thus, the number of disjoint and not $B$--isotopic annuli carried by $N(B)\cap N$ is bounded by a number which depends only on the number of components of $\partial_hN(B)\cap N$.  So, by Scharlemann's theorem and the arguments above, for any $S_n\cap N$, there is a finite collection of components of $S_n\cap N$, denoted by $A_1,\dots, A_k$, such that each component of $S_n\cap N$ is $B$--isotopic to some $A_i$ and $k$ is bounded from above by a fixed number depending only on $N(B)\cap N$.  So, we can split $B$ is a neighborhood of $N$ so that, after the splittings, $B\cap N$ becomes a collection of disjoint compact surfaces $A_i$'s above.  By the assumptions on the $A_i$'s, the branched surface after this splitting still carries $S_n$ and clearly does not carry $\mu$, since $N(B_\mu)\subset N$.  

Suppose $\{S_n\}$ is the sequence of strongly irreducible Heegaard surfaces fully carried by $B$.  Then, for each $S_n$, we use $\Sigma_n$ to denote the union of those $A_i$'s above.  Now, we consider the sequence $\{\Sigma_n\}$.  Each $\Sigma_n$ is fully carried by $B\cap N$ and consists of $\partial$--parallel annuli plus at most one exceptional component.  Since the number of components of $\Sigma_n$ is bounded by a fixed number and since there are only finitely many equivalence classes,  $\{\Sigma_n\}$ belong to finitely many isotopy classes.  So, if we split $B\cap N$ into the $A_i$'s for each $n$, we only get a finite number of different branched surfaces, up to isotopy.  Therefore, we can split $B$ in a neighborhood of $N$ and obtain a finite collection of branched surfaces such that
\begin{enumerate}
\item up to isotopy in $N$, each $S_n$ is fully carried by a branched surface in this collection,
\item the intersection of $N$ and each branched surface in this collection consists of annuli and at most one exceptional component as in Scharlemann's theorem.  In particular, no branched surface in this collection carries $\mu$.
\end{enumerate}

By combining all the splittings before, we get a finite collection of branched surfaces satisfying the properties of Lemma~\ref{Lmu}.
\end{proof}

Using Lemma~\ref{Ltorus2} and Theorem~\ref{TMS}, we can drop the hypothesis that $\mu$ is an exceptional minimal lamination in Lemma~\ref{Lmu}.

\begin{corollary}\label{Cmu}
Let $B$, $B_N$ and $\{S_n\}$ be as in Lemma~\ref{Lmu}.  Let $\mu$ be a measured lamination carried by $B_N$ with $\chi(\mu)=0$.  Then, $B$ can be split into a finite collection of branched surfaces with the following properties.
\begin{enumerate}
\item Up to isotopy, each $S_n$ is carried by a branched surface in this collection.
\item No branched surface in this collection carries $\mu$.
\end{enumerate}
\end{corollary}
\begin{proof}
By Theorem~\ref{TMS}, $\mu$ is a disjoint union of a finite number of sub-laminations, $\mu_1,\dots,\mu_m$.  It is a well-known fact that a measured lamination with positive Euler characteristic has a 2--sphere (or $P^2$) leaf (this is even true for ``abstract" laminations, see \cite{MO}).  Since $B_N$ does not carry any 2--sphere, $B_N$ does not carry any measured lamination with positive Euler characteristic.  So, $\chi(\mu_i)=0$ for each $i$.  By Corollary~\ref{Cklein}, $B_N$ does not carry any Klein bottle.  Hence, each $\mu_i$ either is an exceptional minimal lamination or consists of parallel normal tori.  Now, the corollary follows immediately from Lemmas~\ref{Ltorus}, \ref{Ltorus2} and \ref{Lmu}.
\end{proof}

\section{Proof of the main theorem}\label{Sproof}

Let $B'$ be a branched surface obtained by splitting $B$.  By section~\ref{Spre}, we may naturally assume $N(B')\subset N(B)$.  Recall that we say that a lamination $\mu$  carried by $B$ is also carried by $B'$ if after some $B$--isotopy, $\mu\subset N(B')\subset N(B)$, transverse to the $I$--fibers.

\begin{proposition}\label{neighbor}
Let $B'$ be a branched surface obtained by splitting $B$.  Suppose $\mu$ is a measured lamination carried by $B$ but not carried by $B'$.  Then, there is a neighborhood $N_\mu$ of $\mu$ in the projective lamination space of $B$, such that no measured lamination in $N_\mu$ is carried by $B'$.
\end{proposition}
\begin{proof}
Let $\mathcal{PL}(B)$ and $\mathcal{PL}(B')$ be the projective lamination spaces for $B$ and $B'$ respectively. Suppose there is a measured lamination carried by $B'$ in every neighborhood of $\mu$ in $\mathcal{PL}(B)$.  Then, there are an infinite sequence of measured laminations $\{\mu_n\}$ carried by $B'$ and the limit point of $\{\mu_n\}$ in $\mathcal{PL}(B)$ is $\mu$.  Since $\mathcal{PL}(B')$ is compact, this sequence $\{\mu_n\}$ must have an accumulation point $\mu'$ in $\mathcal{PL}(B')$.  Since every lamination carried by $B'$ is also carried by $B$, $\mu'$ is carried by $B$ and hence is an accumulation point of $\{\mu_n\}$ in $\mathcal{PL}(B)$.  So, $\mu$ and $\mu'$ correspond to the same point in $\mathcal{PL}(B)$.  Since points in $\mathcal{PL}(B)$ and the measured laminations described in section~\ref{SPL} are one-to-one correspondent, $\mu=\mu'$ and we get a contradiction.
\end{proof}

Now, we are in position to prove Theorem~\ref{T2}.

\begin{proof}[Proof of Theorem~\ref{T2}]
By \cite{BCZ,BO,M,MSch}, we may assume $M$ is not a lens space or a small Seifert fiber space.  So, by \cite{JR}, $M$ admits a 0--efficient triangulation.   By \cite{R, St, K}, every strongly irreducible Heegaard splitting can be isotoped to an almost normal surface.  By section~\ref{Spre}, we can find a finite collection of branched surfaces, $B_1,\dots,B_n$, such that each almost normal Heegaard surface is fully carried by some $B_i$.  Since the triangulation is 0--efficient, we may assume that $B_i$ does not carry any normal 2--sphere for each $i$.

Let $\mathcal{PL}(B_i)$ be the projective lamination space for $B_i$.  We can identify each point in $\mathcal{PL}(B_i)$ to a measured lamination carried by $B_i$.  Let $\mathcal{T}_i\subset\mathcal{PL}(B_i)$ be the collection of normal measured laminations with Euler characteristic 0.  By Proposition~\ref{Pcompact}, $\mathcal{T}_i$ is a compact subset of $\mathcal{PL}(B_i)$.  

By Corollary~\ref{Cmu}, for each normal measured lamination $\mu\in\mathcal{T}_i$ carried by $B_i$, we can split $B_i$ into a finite collection of branched surfaces such that any strongly irreducible Heegaard surface fully carried by $B_i$ is fully carried by a branched surface in this collection, and no branched surface in this collection carries $\mu$.  Since this is a finite collection, by Proposition~\ref{neighbor}, there is a neighborhood $N_\mu$ of $\mu$ in $\mathcal{PL}(B_i)$ such that none of the measured laminations in $N_\mu$ is carried by any branched surface in this collection. 

Since $\mathcal{PL}(B_i)$ and $\mathcal{T}_i$ are compact, there are a finite number of normal measured laminations $\mu_1,\dots,\mu_k$ in $\mathcal{T}_i$ such that $\cup_{j=1}^kN_{\mu_j}$ covers $\mathcal{T}_i$.  By applying Corollary~\ref{Cmu} to each $\mu_j$ and combining all the splittings for the $\mu_j$'s, we can split $B_i$ into a finite collection of branched surfaces such that 
\begin{enumerate}
\item each strongly irreducible Heegaard surface fully carried by $B_i$ is still fully carried by a branched surface in this collection,
\item no branched surface in this collection carries any measured lamination in $N_{\mu_j}$ ($j=1,\dots, k$).  Since $\cup_{j=1}^kN_{\mu_j}$ covers $\mathcal{T}_i$, no branched surface in this collection carries any normal measured lamination with Euler characteristic 0.
\end{enumerate}

After performing such splitting on each $B_i$, we get a finite collection of branched surfaces with the desired properties.
\end{proof}

Now, Theorem~\ref{main} follows easily from Theorem~\ref{T2}.

\begin{proof}[Proof of Theorem~\ref{main}]
By Theorem~\ref{T2}, there are finitely many branched surfaces, say $B_1,\dots,B_n$, such that any almost normal strongly irreducible Heegaard surface is fully carried by some $B_i$, and for any $i$, $B_i$ does not carry any normal 2--sphere or normal torus.  

Since each almost normal surface has at most one almost normal piece, at most one branch sector of $B_i$ contains an almost normal piece, and (if it exists) we call this branch sector the \emph{almost normal sector} of $B_i$.  For each almost normal surface $S$ fully carried by $B_i$, the weight of $S$ at the almost normal sector is exactly one.  It is possible that $B_i$ has no almost normal sector, in which case every surface carried by $B_i$ is normal.

Let $\mathcal{S}_i$ be the set of almost normal Heegaard surfaces fully carried by $B_i$ and with a fixed genus $g$.  Each $S\in\mathcal{S}_i$ corresponds to a positive integer solution to the branch equations.  We can write $S=(x_1,\dots,x_m)$, where each $x_i$ is the weight of $S$ at a branch sector of $B_i$.  We may assume the first coordinate $x_1$ corresponds to the almost normal sector.  So, $x_1=1$ for every $S\in\mathcal{S}_i$.

If $\mathcal{S}_i$ is an infinite set, then we can find two surfaces $S_1=(x_1,\dots,x_m)$ and $S_2=(y_1,\dots,y_m)$ in $\mathcal{S}_i$ such that $x_i\le y_i$ for each $i$.  There are many ways to see this and the following is suggested by the referee.  If for a fixed $i$ infinitely many surfaces in $\mathcal{S}_i$ have the same $i$-th coordinate, then we work with that collection and proceed by induction. Eventually we reach an infinite collection of surfaces where for some coordinates they all agree and for the rest only finitely many take any one value. Then for a fixed surface $S_1=(x_1,\dots,x_m)$ there are only finitely many surfaces in this collection with any coordinate less than $\max_i\{x_i\}$. Hence such a surface $S_2=(y_1,\dots,y_m)$ exists.

Thus, $T=(y_1-x_1, y_2-x_2,\dots, y_m-x_m)$ is a non-negative integer solution to the branch equations.  So, we may consider $T$ as a closed surface carried by $B_i$ ($T$ may not be connected).  If $B_i$ has an almost normal sector, then $x_1=y_1=1$ and the weight of $T$ at the almost normal sector is $0$.  Hence, $T$ is a normal surface.  Moreover, since $genus(S_1)=genus(S_2)=g$, $\chi(T)=\chi(S_1)-\chi(S_2)=0$.  This is impossible, since $B_i$ does not carry any normal 2--sphere or normal torus.  Hence, each $\mathcal{S}_i$ is a finite set and there are only finitely many strongly irreducible Heegaard splittings of any genus $g$.  

Johannson proved Theorem~\ref{main} for Haken manifolds \cite{Jo1,Jo2}.  For non-Haken manifolds, by \cite{CG}, a weakly reducible Heegaard splitting is in fact reducible.  So, any weakly reducible Heegaard splitting in an irreducible non-Haken manifold can be destabilized into a strongly irreducible Heegaard splitting.  Therefore, Theorem~\ref{main} holds for all Heegaard splittings
\end{proof}

\end{psfrags}

\end{document}